\input epsf.sty
\input rotate.tex
\input arrow.tex
\baselineskip=14pt

\hsize=118mm 
\hoffset=20mm
\vsize=215mm
\voffset=15mm

\overfullrule=0pt

\font\tencsc=cmcsc10

\font\tenbb=msbm10
\font\sevenbb=msbm7
\font\fivebb=msbm5
\newfam\bbfam
\textfont\bbfam=\tenbb \scriptfont\bbfam=\sevenbb
\scriptscriptfont\bbfam=\fivebb
\def\bb{\fam\bbfam}

\def\Cb{{\bb C}}
\def\Zb{{\bb Z}}
\def\Rb{{\bb R}}

\def\rot#1#2{\rotdimen=\ht#1\advance\rotdimen by\dp#1%
   \hbox to\rotdimen{\hskip\ht#1\vbox to\wd#1{\rotstart{90 rotate}%
   \box#1\vss}\hss}\rotfinish}

\def\sevafig#1#2{\centerline{
 \epsfxsize=#2\epsfbox{#1}}}

\def\build#1_#2^#3{\mathrel{
\mathop{\kern 0pt#1}\limits_{#2}^{#3}}}

\def\diagram#1{\def\normalbaselines{\baselineskip=0pt
\lineskip=10pt\lineskiplimit=1pt}   \matrix{#1}}

\def\hfl#1#2{\smash{\mathop{\hbox to 6mm{\rightarrowfill}}
\limits^{\scriptstyle#1}_{\scriptstyle#2}}}

\def\hfll#1#2{\smash{\mathop{\hbox to 6mm{\leftarrowfill}}
\limits^{\scriptstyle#1}_{\scriptstyle#2}}}

\def\vfll#1#2{\llap{$\scriptstyle #1$}\left\uparrow
\vbox to 3mm{}\right.\rlap{$\scriptstyle #2$}}

\def\Ac{{\cal A}}
\def\Bc{{\cal B}}
\def\Dc{{\cal D}}
\def\Uc{{\cal U}}
\def\Cc{{\cal C}}
\def\Hc{{\cal H}}

\def\ify{\infty}
\def\ra{\rightarrow}
\def\longra{\longrightarrow}
\def\ot{\otimes}
\def\op{\oplus}
\def\bu{\bullet}
\def\wdg{\wedge}
\def\wh{\widehat}
\def\wt{\widetilde}
\def\sm{\simeq}
\def\ov{\overline}

\def\mpo{\mapsto}
\def\sbsq{\subseteq}
\def\bsh{\backslash}
\def\part{\partial}
\def\ts{\times}
\def\sbs{\subset}
\def\lgl{\langle}
\def\rgl{\rangle}

\def\vp{\varphi}
\def\g{\gamma}
\def\s{\sigma}
\def\om{\omega}
\def\a{\alpha}
\def\t{\theta}
\def\b{\beta}

\def\D{\Delta}
\def\Lb{\Lambda}
\def\Om{\Omega}
\def\Si{\Sigma}
\def\G{\Gamma}

\def\xx{\vrule height 0.5em depth 0.2em width 0.5em}


\centerline{\tencsc On the cyclic Formality conjecture}

\bigskip

\centerline{Boris Shoikhet}

\vglue 2cm

\noindent 1. Introduction

\smallskip

\noindent 2. Cyclic Hochschild-Kostant-Rosenberg Theorem

\smallskip

\noindent 3. Cyclic Formality $L_{\ify}$-morphism

\smallskip

\noindent 4. Globalization

\smallskip

\noindent 5. Applications

\vglue 1cm

\noindent {\bf Abstract.} We conjecture an explicit formula for a cyclic 
analog of the Formality $L_{\ify}$-morphism [K]. We prove that its first 
Taylor component, the cyclic Hochschild-Kostant-Rosenberg map, is in fact 
a morphism (and a quasiisomorphism) of the complexes. To prove it we 
construct a cohomological version of the Connes-Tsygan bicomplex in 
cyclic homology. As an application of the cyclic Formality conjecture, we 
obtain an explicit formula for cyclically invariant deformation 
quantization. We show that (a more precise version of) the Connes-Flato-Sternheimer 
conjecture [CFS] on the existence of closed star-products on 
a symplectic manifold also follows from our conjecture.

\vglue 1cm

\noindent {\bf 1. Introduction}

\smallskip

\noindent {\bf 1.1.} Here we recall the definition of the cohomological
Hochschild complex of an associative algebra $A$, and the definitions of 
the associative product on this complex and the Gerstenhaber bracket. Let 
$A$ be an associative algebra over $\Cb$. The Hochschild complex is 
complex
$$
0 \ra C^{0} \build \ra_{}^{d} C^{1} \build \ra_{}^{d} C^{2} \build \ra_{}^{d} 
\ldots
$$
where $C^{k} = {\rm Hom}_{\Cb} (A^{\ot k} , A) \ , \ k \geq 0$, and for 
$\vp \in C^{k}$ the differential $d\vp \in C^{k + 1}$ is defined as 
follows:
$$
\eqalignno{
&(d\vp) (a_1 \ot \ldots \ot a_{k + 1}) = a_1 \cdot \vp (a_2 \ot \ldots 
\ot a_{k + 1}) - \vp ((a_1 \cdot a_2) \ot a_3 \ot \ldots \ot a_n ) \cr
+ & \ \vp (a_1 \ot (a_2 \, a_3) \ot a_4 \ot \ldots \ot a_{n+1}) - \ldots
\pm \vp (a_1 \ot a_2 \ot \ldots \ot (a_n \cdot a_{n+1})) \cr
\mp & \ \vp (a_1 \ot \ldots \ot a_n) \cdot a_{n+1} \, . &(1) \cr
}
$$

There exists a nice interpretation of the last formula as well as of the 
Gerstenhaber bracket. Let $\Ac$ be the cofree coalgebra, cogenerated by 
the vector space $A[1]$. As a vector space, $\Ac = \Cb \op A[1] \op A[1]^{\ot 2} 
\op A[1]^{\ot 3} \op \ldots$, and the comultiplication $\D$ is defined as 
follows:
$$
\eqalignno{
&\D (a_1 \ot \ldots \ot a_{k}) = \cr
&1 \ot (a_1 \ot \ldots \ot a_{k}) + a_1 \ot (a_2 \ot \ldots \ot a_{k}) + 
\ldots &(2) \cr
&\ldots + (a_1 \ot \ldots \ot a_{k - 1}) \ot a_{k} + (a_1 \ot \ldots \ot 
a_{k}) \ot 1 \, . \cr
}
$$
Any map $D : A^{\ot k} \ra A$ defines uniquely a coderivation of the 
coalgebra $\Ac$ (as well as a map $B \ra B^{\ot k}$ defines uniquely the 
derivation of the tensor algebra $\Bc$, generated by the vector space 
$B$). We still denote by $D$ the corresponding coderivation of the  
coalgebra $\Ac$. The bracket $[D_1 , D_2]$ of the two coderivations 
defined by maps $D_1 : A^{\ot k} \ra A$ and $D_2 : A^{\ot \ell} \ra A$  
defines a map $[D_1 , D_2] : A^{\ot (k + \ell - 1)} \ra A$ and the 
last map is called the Gerstenhaber bracket of $D_1$ and $D_2$. Moreover, 
if the space $A$ has an algebra structure, there exists the canonical 
coderivation $m : A^{\ot 2} \ra A$, $m (a_1 \ot a_2) = a_1 \cdot a_2$. 
One can check that $[m,m] = 0$. Therefore, the map $d : C^{k} \ra 
C^{k+1}$,
$$
d(\psi)=[m,\psi]
\eqno(3)
$$ 
 satisfies the differential equation $d^2 = 0$. One can check 
that this is exactly the Hochschild differential given by formula (1). 
One can check that the Gerstenhaber bracket defines a $\Zb$-graded Lie 
algebra structure on the Hochschild complex $C^{\bu}$. Then it follows 
from the definition (3) that $C^{\bu}$ equipped with the differential $d$ 
and the Gerstenhaber bracket defines a $dg$ Lie algebra structure on the 
complex $C^{\bu} [+1]$.

In the sequel we will need the explicit formula for the Gerstenhaber 
bracket. Let $\psi_1 \in C^{k_1 + 1}$, $\psi_2 \in C^{k_2 + 1}$. We have:
$$
[\psi_1 , \psi_2] = \psi_1 \circ \psi_2 - (-1)^{k_1 k_2} \, \psi_2 \circ 
\psi_1 \eqno (4)
$$
where
$$
\eqalignno{
&(\psi_1 \circ \psi_2) (a_0 \ot \ldots \ot a_{k_1 + k_2}) = \cr
= & \ \sum_{i=0}^{k_1} (-1)^{ik_2} \, \psi_1 (a_0 \ot \ldots \ot a_{i-1} 
\ot \psi_2 (a_i , a_{i+1} , \ldots , a_{i+k_2}) \ot \cr
& \ \ot a_{i+k_2 + 1} \ot \ldots \ot a_{k_1 + k_2}) \, . &(5) \cr
}
$$

It is a standard fact that the cohomology of the complex  $C^{\bu}$ are equal to
$$
H^k (C^{\bu}) = {\rm Ext}_{A \ot A^0}^k (A,A) \, , \eqno (6)
$$
the $k$-th Ext's group in the category of $A$-bimodules, where the 
algebra $A$ is considered as an $A$-bimodule. According to (6), there 
exists a canonical associative product on the cohomology $H^{\bu} 
(C^{\bu})$. In fact, this product also can be defined on the complex 
$C^{\bu}$ as follows: for $\vp \in C^{k_1}$, $\psi \in C^{k_2}$
$$
(\vp \cdot \psi) (a_1 \ot \ldots \ot a_{k_1 + k_2}) = \vp (a_1 \ot \ldots 
\ot a_{k_1}) \cdot \psi (a_{k_1 + 1} \ot \ldots \ot a_{k_2}) \, . \eqno 
(7)
$$

The associative product and the Lie bracket on the complex $C^{\bu}$ are 
not compatible, i.e. $[\psi_1 \cdot \psi_2 , \psi_3] \ne [\psi_1 , 
\psi_3] \cdot \psi_2 \pm \psi_1 \cdot [\psi_2 , \psi_3]$.

\bigskip

\noindent {\bf 1.2.} Let $A = \Cb \, [x_1 , \ldots , x_d]$ or $A = 
C^{\ify} (M)$ where $M$ is a smooth manifold. Then the cohomology $H^i 
(C^{\bu})$ (the cohomology of the polydifferential part of $C^{\bu}$ in the 
smooth 
case) 
are equal to the space of polynomial polyvector fields on 
$\Rb^d$ or smooth polyvector fields on the manifold $M$. The induced 
associative product and bracket are exactly the $\Lb$-product of 
polyvector fields and the Schouten-Nijenhuis bracket of polyvector 
fields, defined as follows: for $k , \ell \geq 0$
$$
\eqalignno{
&[\xi_0 \wdg \ldots \wdg \xi_k , \eta_0 \wdg \ldots \wdg \eta_{\ell}] = 
&(8) \cr 
= & \ \sum_{i=0}^{k} \sum_{j=0}^{\ell} (-1)^{i+j+k} [\xi_i , \eta_j] \wdg 
\xi_0 \wdg \ldots \wh{\xi}_i \wdg \ldots \wdg \xi_k \wdg \eta_0 \wdg 
\ldots \wdg \wh{\eta}_j \wdg \ldots \eta_{\ell} \cr
}
$$
where $\{ \xi_i \}$ and $\{ \eta_j \}$ are vector fields.

Denote by $T_{\rm poly}^i$ the space of $(i+1)$-vector fields, and by 
$\Dc_{\rm poly}^i$ the polydifferential part of the space ${\rm 
Hom}_{\Cb} (A^{\ot (i+1)} , A)$. 

The Hochschild-Kostant-Rosenberg map $\vp_{\rm HKR} : T_{\rm poly}^{\bu} 
\ra \Dc_{\rm poly}^{\bu}$ is defined as follows:
$$
\vp_{\rm HKR} (\xi_1 \wdg \ldots \wdg \xi_k) (f_1 \ot \ldots \ot f_k) = 
{1 \over k!} \ \build {\rm Alt}_{\xi_1 , \ldots , \xi_k}^{} \ \xi_1 (f_1) 
\cdot \ldots \cdot \xi_k (f_k) \, .
$$

\medskip

\noindent {\bf Theorem.} (Hochschild-Kostant-Rosenberg) 

\smallskip

\item{(i)} {\it the map $\vp_{\rm HKR} : T_{\rm poly}^{\bu} \ra \Dc_{\rm 
poly}^{\bu}$ is a quasiisomorphism of the complexes;}

\smallskip

\item{(ii)} {\it the induced map $T_{\rm poly}^{\bu} \ra H^{\bu} (\Dc_{\rm 
poly}^{\bu})$ is an (iso)morphism of both associative and Lie algebras.} \hfill 
\xx

\bigskip

It follows from this theorem that for any $\eta_1 , \eta_2 \in T_{\rm 
poly}^{\bu}$ one has
$$
[\vp_{\rm HKR} (\eta_1) , \vp_{\rm HKR} (\eta_2)] = \vp_{\rm HKR} ( [\eta_1 
, 
\eta_2]) \ \hbox{\it mod coboundaries}.
$$

The Formality theorem of M.~Kontsevich [K] states that the $dg$ Lie 
algebras $T_{\rm poly}^{\bu}$ and $\Dc_{\rm poly}^{\bu}$ are 
quasiisomorphic 
$dg$ Lie algebras, i.e. there exists a $dg$ Lie algebra ? and diagram
$$
\matrix{
&&? \cr
&{\scriptstyle s} \nearrow &&\nwarrow {\scriptstyle t} \cr
T_{\rm poly}^{\bu} &&&&\Dc_{\rm poly}^{\bu} \cr
}
$$
where the maps $s$ and $t$ are quasiisomorphisms of the $dg$ Lie algebras. 
In fact, this result was proved using the language of the homotopical algebra, 
and it was constructed an $L_{\ify}$-quasiisomorphism $\Uc : T_{\rm 
poly}^{\bu} \ra \Dc_{\rm poly}^{\bu}$ (see [K]).

The analogous result also holds for the associative $dg$ algebras $T_{\rm 
poly}^{\bu}$ and $\Dc_{\rm poly}^{\bu}$ (see [Sh]).

\bigskip

\noindent {\bf 1.3.} Let $A = \Cb \, [x_1 , \ldots , x_d]$. To formulate the 
cyclic Formality conjecture, we need an additional data -- a volume form 
$\Om$ on the space $\Rb^d$. Let us suppose that this form is fixed.

\medskip

\noindent {\bf 1.3.1.} The form $\Om$  
defines an isomorphism $T_{\rm poly}^{i-1} \ra \Om_{\rm DR}^{d-i}$ (recall, 
that $T_{\rm poly}^{i-1}$ is the space of $i$-polyvector fields).

\medskip

\noindent {\bf Definition.} (The divergention operator.) The map ${\rm div} 
: T_{\rm poly}^i \ra T_{\rm poly}^{i-1}$ is defined as follows: $T_{\rm 
poly}^i \build{\wt{\longra}}_{}^{\Om} \Om_{\rm DR}^{d-i-1} \build 
\longra_{}^{d_{\rm DR}} \Om^{d-i} \build{\wt{\longra}}_{}^{\Om} T_{\rm 
poly}^{i-1}$.

The analog of the $dg$ Lie algebra $T_{\rm poly}^{\bu}$ for the cyclic 
Formality conjecture is $T_{\rm poly}^{\bu} \ot \Cb \, [u]$, where $\deg u=2$, 
equipped with the differential
$$
d_{\rm div} (\g \ot u^k) = {\rm div} (\g) \ot u^{k+1} \eqno (9)
$$
and the bracket
$$
[\g_1 \ot u^{k_1} , \g_2 \ot u^{k_2}] = [\g_1 , \g_2] \ot u^{k_1 + k_2} \, 
. 
\eqno (10)
$$

It follows from the Lemma below that this is actually a $dg$ Lie algebra.

\medskip

\noindent {\bf Lemma.} {\it For $\eta_1 \in T_{\rm poly}$, $\eta_2 \in T_{\rm 
poly}$ one has}
$$
{\rm div} \, [\eta_1 , \eta_2] = [{\rm div} \, \eta_1 , \eta_2] + 
(-1)^{i+1} 
\, [\eta_1 , {\rm div} \, \eta_2] .
$$

\medskip

\noindent {\it Proof.} One can prove that $$[\eta_1 , \eta_2]
=\pm({\rm div} \, (\eta_1 \wdg 
\eta_2) 
-
({\rm div} \, \eta_1) \wdg \eta_2 -(-1)^{i+1} \, \eta_1 \wdg ({\rm div} \, 
\eta_2)).$$ The statement of Lemma is a direct consequence of this formula.

\hfill \xx

\bigskip

\noindent {\bf 1.3.2.} The cyclic analog of the $dg$ Lie algebra $\Dc_{\rm 
poly}^{\bu}$ is defined in a bit more tricky way.

Let us define an operator of the cyclic shift $C : \Dc_{\rm poly}^{\bu} \ra 
\Dc_{\rm poly}^{\bu}$.

\medskip

\noindent {\bf Definition.}
$$
\int_{\Rb^d} \psi (f_1 , \ldots , f_n) \cdot f_{n+1} \cdot \Om = (-1)^n 
\cdot 
\int_{\Rb^d} C (\psi) (f_2 , \ldots , f_{n+1}) \cdot f_1 \cdot \Om \, , 
\eqno 
(11)
$$
where $f_1 , \ldots , f_{n+1}$ are functions with a compact support. The 
operator $C : \Dc_{\rm poly}^{\bu} \ra \Dc_{\rm poly}^{\bu}$ is defined by 
the 
continuity for any functions $f_1 , \ldots , f_{n+1}$.

\medskip

\noindent {\bf Definition.}
$$
[\Dc_{\rm poly}^i]_{\rm cycl} = \{ \psi \in \Dc_{\rm poly}^i \mid C (\psi) 
= 
\psi \} \, .
$$

\medskip

\noindent {\bf Lemma.} {\it $[\Dc_{\rm poly}^{\bu}]_{\rm cycl}$ is a $dg$ Lie 
algebra with respect to the Hochschild differential and the Gerstenhaber 
bracket.}

We prove the statement of the Lemma for the differential, the proof for the 
Gerstenhaber bracket is analogous. We have:
$$
\eqalign{
& \ \int (d_{\rm Hoch} \, \psi) (f_1 , \ldots , f_{n+1}) \cdot f_{n+2} 
\cdot \Om \cr
= & \ \int (f_1 \, \psi (f_2 , \ldots , f_{n+1}) - \psi (f_1 \cdot f_2 , 
f_3 , \ldots , f_{n+1}) + \ldots - \ldots \cr
& \ \pm \psi (f_1 , \ldots , f_n) \, f_{n+1}) \, f_n \cr
= & \ \int [ f_1 \, f_{n+2} \, \psi (f_2 , \ldots , f_{n+1}) - (-1)^n \, C 
(\psi) (f_3 , \ldots , f_{n+1} , f_{n+2}) \, f_1 \, f_2 \cr
& \ + (-1)^n \, C(\psi) (f_2 \, f_3 , f_n , \ldots , f_{n+2}) \cdot f_1 
\ldots \cr
& \ \pm (-1)^n \, C (\psi) (f_2 , \ldots , f_n , f_{n+1} \, f_{n+2}) \, f_1] 
\cdot \Om = \cr
= & \ (-1)^{n+1} \int (d_{\rm Hoch} \, \psi) (f_2 , \ldots , f_{n+2}) \cdot 
f_1 \, , \cr
}
$$
because
$$
C(\psi) = \psi \, .
$$
\hfill \xx

\bigskip

\noindent {\bf 1.4. Cyclic Formality conjecture} 

\smallskip

\noindent {\bf Conjecture.} {\it The $dg$ Lie algebras $T_{\rm poly}^{\bu} 
\ot \Cb \, [u]$ and $[\Dc_{\rm poly}^{\bu}]_{\rm cycl}$ are quasiisomorphic 
for 
any volume form $\Om$.}

\medskip

In this form the Conjecture is due to M.~Kontsevich (private 
communication). 
In the present paper we construct explicitly set of maps:
$$
\eqalign{
&\Cc_1 : T_{\rm poly}^{\bu} \ot \Cb \, [u] \ra [\Dc_{\rm poly}^{\bu}]_{\rm 
cycl} \cr
&\Cc_2 : \Lb^2 (T_{\rm poly}^{\bu} \ot \Cb \, [u]) \ra [\Dc_{\rm 
poly}^{\bu}]_{\rm cycl} \, [-1] \cr
&\Cc_3 : \Lb^3 (T_{\rm poly}^{\bu} \ot \Cb \, [u]) \ra [\Dc_{\rm 
poly}^{\bu}]_{\rm cycl} \, [-2] \cr
}
$$
and conjecture that these maps are the components of an $L_{\ify}$-morphism 
(see (34) below), and also an $L_{\ify}$-quasiisomorphism.

The map $\Cc_1$ is a cyclic analog of the Hochschild-Kostant-Rosenberg map 
$\vp_{\rm HKR}$ (see Section~1.2). Let us note, that even this map is quite 
nontrivial. We prove that the map $\Cc_1$ is in fact a (quasiiso)morphism 
of the 
complexes.

Let us note that the map $C : \Dc_{\rm poly}^i \ra \Dc_{\rm poly}^i$ 
satisfies the equality
$$
C^{i+2} = 1 \, . \eqno (12)
$$
We denote
$$
\sum = 1+C+\ldots + C^{i+1} \, . \eqno (13)
$$
It is clear that 
$$
C \left( \sum \psi \right) = \sum  \psi \eqno (14)
$$
for any $\psi$. The map $\sum$ is not compatible with the differential and 
the Gerstenhaber bracket, but, in a sense, it almost is. We will 
discuss the compatibility with the Hochschild differential in the next Section.
Compatibility, in a certain sense, of the map $\sum$ with the Gerstenhaber 
bracket 
is the most mysterious part
of the cyclic Formality conjecture, and the main Conjecture 3.2.2
 can be considered as an expression of this
compatibility.

\vglue 1cm

\noindent {\bf 2. Cyclic Hochschild-Kostant-Rosenberg Theorem}

\smallskip

\noindent {\bf 2.1.}  We consider the $dg$ Lie algebras
$$
\{ T_{\rm poly}^{\bu} \ot \Cb \, [u] , d_{\rm div} \} \quad \hbox{and} \quad 
\{ 
[\Dc_{\rm poly}^{\bu}]_{\rm cycl} , d_{\rm Hoch} \}
$$
constructed from the algebra $A = \Cb \, [x_1 ,\ldots ,x_d]$ (see Section 
1.3). 
We denote $[T_{\rm poly}^{\bu}]_{\rm div} = \{ \g \in T_{\rm poly}^{\bu} 
\mid {\rm div} \, \g = 0 \}$. It follows from the Poincar\'e lemma that
$$
H^i (T_{\rm poly}^{\bu} \ot \Cb \, [u]) \sm \left\{ \matrix{
[T_{\rm poly}^i]_{\rm div}, &i \leq d-1 \hfill \cr
\Cb , \hfill &i=d+2k-1 , &k \geq 0 \, . \cr
} \right. \eqno (15)
$$

The main result of this Section is the following

\medskip

\noindent {\bf Theorem.} 
$$
H^i ([\Dc_{\rm poly}^{\bu}]_{\rm cycl}) \sm \left\{ \matrix{
[T_{\rm poly}^i]_{\rm div}, &i \leq d-1 \hfill \cr
\Cb , \hfill &i=d+2k-1 , &k \geq 0 \, . \cr
} \right.
$$

\medskip

\noindent {\bf 2.1.1.} The analogous result in cyclic homology was proved 
using the spectral sequence, connected with the Connes-Tsygan bicomplex. It 
turns out that there exists an explicit analog of this construction in our 
situation.

Let $A$ be an associative algebra. We define the complex $\{ K^{\bu} , d_K 
\}$ as follows:
$$
K^i = {\rm Hom}_{\Cb} (A^{\ot i} , A) \ , \ i \geq 0 \, , \eqno (16)
$$
and for $\psi \in K^i$
$$
\eqalignno{
& \ (d_K \, \psi) (a_1 \ot \ldots \ot a_{i+1}) = a_1 \, \psi (a_2 \ot 
\ldots 
\ot a_{i+1}) - \cr
- & \ \psi (a_1 a_2 \ot a_3 \ot \ldots \ot a_{i+1}) + \ldots \pm \psi (a_1 
\ot \ldots \ot a_i a_{i+1}) \, . &(17) \cr
}
$$

Note that it is exactly the Hochschild differential without the last term.

\medskip

\noindent {\bf Lemma.}

\smallskip

 {\it $d_K^2 : K^i \ra K^{i+2}$, $i \geq 0$, is equal to $0$ for 
any associative algebra $A$.}

\medskip

\noindent {\it Proof.}  Is straightforward. \hfill \xx

\bigskip

\noindent {\bf 2.1.2. Lemma.} {\it For any associative algebra $A$ with 
unit 
the complex $\{ K^{\bu} , d_K \}$ is acyclic.}

\medskip

\noindent {\it Proof.} It is sufficient to construct a homotopy $h : 
K^{\bu} \ra K^{\bu - 1}$ such that
$$
d_K \, h \pm h \, d_K = \pm \, {\rm Id} \, . \eqno (18)
$$
We set
$$
(h \, \psi) (a_1 \ot \ldots \ot a_{k-1}) = \psi (a_1 \ot \ldots \ot a_{k-1} 
\ot 1) \, . \eqno (19)
$$
It is clear that (18) is satisfied. \hfill \xx

\bigskip

\noindent {\bf 2.1.3.} Let $A$ be the algebra of functions on a smooth 
manifold $M$, $A = C^{\ify} (M)$, and $\Om$ be a volume form on $M$. Then 
formula (11) defines the operator $C : \Dc_{\rm poly}^i (M) \ra \Dc_{\rm 
poly}^i (M)$ such that $C^{i+2} = 1$, and we set $\sum = 1+C+C^2 + \ldots + 
C^{i+1}$.

The results of Sections 2.1.1, 2.1.2 are true also for polydifferential 
part of the complex $\{ K^{\bu} , d_K \}$. We denote by ${\rm Hom}_{\Cb}^{\rm 
poly} (A^{\ot k} , A)$ the polydifferential part of ${\rm Hom}_{\Cb} (A^{\ot k} 
, A)$.

\medskip

\noindent {\bf Key-lemma.} {\it The following is a bicomplex:}
$$
\diagram{
&&\vfll{}{} &&\vfll{}{} &&\vfll{}{} \cr
0 &\hfl{}{} &{\rm Hom}_{\Cb}^{\rm poly} (A^{\ot 3} , A) &\hfl{1-C}{} 
&{\rm Hom}_{\Cb}^{\rm poly} (A^{\ot 3} , A) &\hfl{\sum}{} &{\rm 
Hom}_{\Cb}^{\rm poly} (A^{\ot 3} , A) &\hfl{1-C}{} \cr
&&\vfll{}{} &\star_1 &\vfll{}{d_K} &\star_2 &\vfll{}{d_{\rm Hoch}} \cr
0 &\hfl{}{} &{\rm Hom}_{\Cb}^{\rm poly} (A^{\ot 2} , A) &\hfl{1-C}{} &{\rm 
Hom}_{\Cb}^{\rm poly} (A^{\ot 2} , A) &\hfl{\sum}{} &{\rm Hom}_{\Cb}^{\rm 
poly} (A^{\ot 2} , A) &\hfl{1-C}{} \cr
&&\vfll{}{d_{\rm Hoch}} &&\vfll{}{d_K} &&\vfll{}{d_{\rm Hoch}} \cr
0 &\hfl{}{} &{\rm Hom}_{\Cb}^{\rm poly} (A,A) &\hfl{1-C}{} &{\rm 
Hom}_{\Cb}^{\rm poly} (A,A) &\hfl{\sum}{} &{\rm Hom}_{\Cb}^{\rm poly} (A,A) 
&\hfl{1-C}{} \cr
&&\vfll{}{d_{\rm Hoch}} &&\vfll{}{d_K} &&\vfll{}{d_{\rm Hoch}} \cr
0 &\hfl{}{} &A &\hfl{1-C}{} &A &\hfl{\sum}{} &A &\hfl{1-C}{} \cr
&&\vfll{}{} &&\vfll{}{} &&\vfll{}{} \cr
&&0 &&0 &&0
}
$$

\noindent {\it Proof.} It follows from (12) that $(1-C) \sum = \sum (1-C) = 
0$. The columns are complexes by the definitions. It remains to prove the 
commutativity of squares.

\smallskip

\noindent (i) Let us prove the commutativity of square $\star_1$. For $\psi 
\in {\rm Hom}_{\Cb}^{\rm poly} (A^{\ot k} , A)$ one has:
$$
\eqalignno{
& \ \int f_{k+2} \cdot (1-C) \, d_{\rm Hoch} \, \psi \cdot \Om = \cr
= & \ \int f_{k+2} \cdot (1-C) \cdot \{ f_1 \cdot \psi (f_2 \ot \ldots \ot 
f_{k+1}) - \psi (f_1 f_2 \ot f_3 \ot \ldots) \cr
+ & \ \ldots \pm \psi (f_1 \ot f_2 \ot \ldots \ot f_k f_{k+1}) \mp \psi 
(f_1 \ot  \ldots \ot f_k) f_{k+1} \} \, \Om \, . &(20) \cr
}
$$
The right-hand side of (20) is equal to
$$
\eqalignno{
& \ \int f_{k+2} \cdot \{ f_1 \psi (f_2 \ot \ldots \ot f_{k+1}) - \psi (f_1 
f_2 \ot \ldots \ot f_{k+1}) + \ldots \cr
\pm & \ \psi (f_1 \ot f_2 \ot \ldots \ot f_k f_{k+1}) \mp \underline{\psi 
(f_1 \ot \ldots \ot f_k) \, f_{k+1}} \} \, \Om \cr
- & \ (-1)^{k-1} \int \{ \underline{f_{k+1} f_{k+2} \psi (f_1 \ot \ldots 
\ot 
f_k \ot f_{k+1})} \cr
- & \ f_{k+1} \, \psi (f_{k+2} \, f_1 \ot f_2 \ot \ldots \ot f_k) + \ldots 
\cr
\pm & \ f_{k+1} \, \psi (f_{k+2} \ot f_1 \ot \ldots \ot f_{k-1} f_k) \cr
\mp & \ f_{k+1} \, \psi (f_{k+2} \ot f_1 \ot \ldots \ot f_{k-1}) \ot f_k \} 
\, \Om \, . &(21)
}
$$

Let us note that underlined terms in (21) cancel each other.

On the other hand,
$$
\eqalignno{
& \ \int f_{k+2} \cdot d_K (1-C) \, \psi \cdot \Om = \cr
= & \ \int f_{k+2} (d_K \, \psi - d_K \, C \psi) \cdot \Om = \cr
= & \ \int f_{k+2} \{ f_1 \, \psi (f_2 \ot \ldots \ot f_{k+1}) - \psi (f_1 
f_2 \ot \ldots \ot f_{k+1}) + \ldots \cr
\pm & \ \psi (f_1 \ot f_2 \ot \ldots \ot f_k f_{k+1}) \} \cdot \Om \cr
- & \ \int f_{k+2} \{ f_1 \cdot (C \psi) (f_2 \ot \ldots \ot f_{k+1}) - (C 
\psi) (f_1 f_2 \ot \ldots \ot f_k) \cr
+ & \ \ldots \pm (C \psi) (f_1 \ot f_2 \ot \ldots \ot f_k f_{k+1}) \} \cdot 
\Om \, . &(22) \cr
}
$$
The second summand in the r.h.s. of (22) is equal to
$$
\eqalign{
- & \ (-1)^k \int \{ f_{k+1} \, \psi (f_{k+2} f_1 \ot \ldots \ot f_k) - 
f_{k+1} \, \psi (f_{k+2} \ot f_1 f_2 \ot \ldots \ot f_k) \cr
& \ + \ldots \pm f_k f_{k+1} \, \psi (f_{k+2} \ot f_1 \ot f_2 \ot \ldots 
\ot 
f_{k-1}) \} \, \Om \, . \cr
}
$$
We see that (21) $=$ (22), and the commutativity of square $\star_1$ is 
proved;

\smallskip

\noindent (ii) the proof of the commutativity of the square $\star_2$ is 
analogous. \hfill \xx

\smallskip

\noindent {\it Remark.} The results of this Subsection 
hold for any associative algebra $A$ with unit equipped with a 
trace functional $\int : A \ra \Cb$, i.e. $\int a \cdot b = \int b \cdot a$ 
for any $a,b \in A$, provided the condition that for any cochain $\psi (a_1 
\ot \ldots \ot a_k)$ there exists the unique cochain $(C\psi) (a_1 \ot 
\ldots \ot a_k)$ such that
$$
\int \psi (a_1 \ot \ldots \ot a_n) \cdot a_{n+1} = (-1)^n \int (C\psi) (a_2 
\ot \ldots \ot a_{n+1}) \cdot a_1 \, .
$$

\bigskip

\noindent {\bf 2.1.4.} Here we prove the following statement.

\medskip

\noindent {\bf Theorem.} {\it Let $M$ be a smooth manifold.}
$$
H^i ([\Dc_{\rm poly}^{\bu} (M)]_{\rm cycl}) \sm \left\{ \matrix{
[T_{\rm poly}^i (M)]_{\rm div} \op H_{\rm DR}^{d-i+1} (M) \op H_{\rm DR}^{d-i+3} 
(M) \op \ldots, \ i \leq d-1 \hfill \cr
\cr
H_{\rm DR}^{\rm even} (M) , \quad i = d-1+k, \ k \ \hbox{even} \hfill \cr
\cr
H_{\rm DR}^{\rm odd} (M) , \quad i = d-1+k, \ k \ \hbox{odd} \, . \hfill \cr
} \right.
$$

\smallskip

\noindent {\it Proof.} The rows of the bicomplex 2.1.3 are acyclic, 
except 
degree 0, because its cohomology are equal to the group cohomology $H^{\bu} 
(\Zb / (n+1) \Zb$, ${\rm Hom}_{\Cb}^{\rm poly} (A^{\ot n} , A))$, and it is 
clear that the last cohomology is zero except $H^0 (\Zb / (n+1) \Zb , {\rm 
Hom}_{\Cb}^{\rm poly} (A^{\ot n} , A)) \sm [{\rm Hom}_{\Cb}^{\rm poly} 
(A^{\ot n} , A)]_{\rm cycl} = [\Dc_{\rm poly}^{n-1} (M)]_{\rm cycl}$. 
Therefore, the bicomplex 2.1.3 is quasiisomorphic to the complex $[\Dc_{\rm 
poly}^{\bu}$ \break $(M)]_{\rm cycl}$.

On the other hand, the second filtration of the bicomplex 2.1.3 gives us 
the 
spectral sequence with second term
$$
\varrowlength=20pt
\commdiag{
\vdots&\vdots&\vdots&\vdots&\vdots&\cr
T_{\rm poly}^2 &0 &T_{\rm poly}^2 &0 &
T_{\rm poly}^2&\ \cdots \cr
&\arrow(3,-1)\rt{d_2}&&\arrow(3,-1)\rt{d_2}\cr
T_{\rm poly}^1 &0 &T_{\rm poly}^1 &0 &
T_{\rm poly}^1&\ \cdots \cr
&\arrow(3,-1)\rt{d_2}&&\arrow(3,-1)\rt{d_2}\cr
T_{\rm poly}^0 &0 &T_{\rm poly}^0 &0 &
T_{\rm poly}^0&\ \cdots \cr
&\arrow(3,-1)\rt{d_2}&&\arrow(3,-1)\rt{d_2}\cr
T_{\rm poly}^{-1} &0 &T_{\rm poly}^{-1} &0 &
T_{\rm poly}^{-1}&\ \cdots \cr
&\arrow(3,-1)&&\arrow(3,-1)\cr
&&0&&0}
$$
because of Lemma 2.1.2.

Here $d_2 = {\rm div} : T_{\rm poly}^i \ra T_{\rm poly}^{i-1}$. 
One can prove that the spectral sequence collapses in the second term. Therefore, 
we obtain the statement of the Theorem. \hfill \xx

\bigskip

\noindent {\bf 2.2. Cyclic Hochschild-Kostant-Rosenberg map (I)}

\smallskip

We want to construct a map $\vp_{\rm HKR}^{\rm cycl} : \{ T_{\rm 
poly}^{\bu} 
\ot \Cb \, [u] , d_{\rm div} \} \ra \{ [\Dc_{\rm poly}^{\bu}]_{\rm cycl}$, 
$d_{\rm Hoch} \}$ which is a quasiisomorphism of the complexes, in the case 
$A = \Cb \, [ x_1 ,$ $\ldots , x_d]$. Here we consider first examples.

\smallskip

\item{(i)} for $f \in T_{\rm poly}^{-1}$ we set
$$
\vp_{\rm HKR}^{\rm cycl} (f) = f \in \Dc_{\rm poly}^{-1} = [\Dc_{\rm 
poly}^{-1}]_{\rm cycl}
$$

\item{(ii)} for $\g = \xi_1 \wdg \ldots \wdg \xi_k \in T_{\rm poly}^{k-1}$ 
we want to define $\vp_{\rm HKR}^{\rm cycl} (\g)$ such that for $\g \in 
[T_{\rm poly}^{k-1}]_{\rm div}$ one has
$$
\vp_{\rm HKR}^{\rm cycl} (\g) = {1 \over k!} \ \build{\rm Alt}_{\xi_1 , 
\ldots , \xi_k}^{} \ \xi_1 (f_1) \cdot \ldots \cdot \xi_k (f_k) = \vp_{\rm 
HKR} (\g) \, .
$$
Let us consider the first case, $\g = \xi \in T_{\rm poly}^0$. We set:
$$
\vp_{\rm HKR}^{\rm cycl} (\xi) (f) = \xi (f) + {1 \over 2} \, {\rm div} 
(\xi) \cdot f = {1 \over 2} \, (\vp_{\rm HKR} (\xi) + C (\vp_{\rm HKR}) 
(\xi)) (f) \, .
$$
It is easy to see that
$$
\int f_1 \cdot \vp_{\rm HKR}^{\rm cycl} (f_2) \cdot \Om = - \int f_2 \cdot 
\vp_{\rm HKR}^{\rm cycl} (f_1) \cdot \Om \, .
$$

\smallskip

\item{(iii)} We have to define $\vp_{\rm HKR}^{\rm cycl} (f \ot u)$ such 
that
$$
\vp_{\rm HKR}^{\rm cycl} (d_{\rm div} (\xi)) = d_{\rm Hoch} (\vp_{\rm 
HKR}^{\rm cycl} (\xi)) \, . \eqno (23)
$$
We set:
$$
\vp_{\rm HKR}^{\rm cycl} (f \ot u) (f_1 \ot f_2) = {1 \over 2} \, f \cdot 
f_1 \cdot f_2 \, .
$$
It is clear that (23) is true.

\smallskip

\item{(iv)} $\g = \xi_1 \wdg \xi_2 \in T_{\rm poly}^1$ we set:
$$
\vp_{\rm HKR}^{\rm cycl} (\g) = {1 \over 3} \, \sum (\vp_{\rm HKR} (\g)) \, 
.
$$
We have:
$$
\eqalign{
& \ \vp_{\rm HKR}^{\rm cycl} (\g) (f_1 \ot f_2) = \vp_{\rm HKR} (\g) (f_1 
\ot f_2) - \cr
- & \ {1 \over 6} \, \big({\rm div} (\xi_1 \wdg \xi_2) (f_1) \cdot f_2 + f_1 
\cdot {\rm div} (\xi_1 \wdg \xi_2) (f_2)\big) \, . \cr
}
$$

\item{(v)}  We want to define
$$
\vp_{\rm HKR}^{\rm cycl} (\xi \ot u) , \ \hbox{where} \quad \xi \in T_{\rm 
poly}^0 \, ,
$$
such that
$$
\vp_{\rm HKR}^{\rm cycl} (d_{\rm div} (\xi_1 \wdg \xi_2)) = d_{\rm Hoch} 
(\vp_{\rm HKR}^{\rm cycl} (\xi_1 \wdg \xi_2)) \, . \eqno (24)
$$
We set:
$$
\eqalign{
& \ \vp_{\rm HKR}^{\rm cycl} (\xi \ot u) \, (f_1 \ot f_2 \ot f_3) = \cr
= & \ {1 \over 6} \, (\xi (f_1) \cdot f_2 \cdot f_3 + f_1 \cdot f_2 \cdot 
\xi (f_3)) + {1 \over 12} \, {\rm div} (\xi) \cdot f_1 \cdot 
f_2 \cdot f_3 \, . \cr
}
$$

\medskip

\noindent {\bf 2.3. Cyclic Hochschild-Kostant-Rosenberg map (II)}

\smallskip

Here we define $\vp_{\rm HKR}^{\rm cycl} (\g \ot u^k)$ for arbitrary $\g 
\ot 
u^k \in T_{\rm poly}^{\bu} \ot \Cb \, [u]$. It is very convenient to use the 
language of graphs of M.~Kontsevich [K].

Let $\G (\ell , k)$ be the set of all the graphs with $\ell + 2k$ vertices 
on the line $\Rb$ (vertices of the second type in [K]), and the unique additional  
vertex (of the first type) and $\ell$ oriented edges started at this vertex 
such that:
$$
\matrix{
\hbox{\it between any two consecutive endpoints} \cr
\hbox{\it of edges there is an even number of other} \cr
\hbox{\it vertices of the second type} \cr
} \eqno (25)
$$
( a vertex of the second type is the endpoint for not more than 1 
edge).

\medskip

\noindent {\it Example.} The set $\G (2,1)$ is shown on the Figure~1

\vglue 1cm

\sevafig{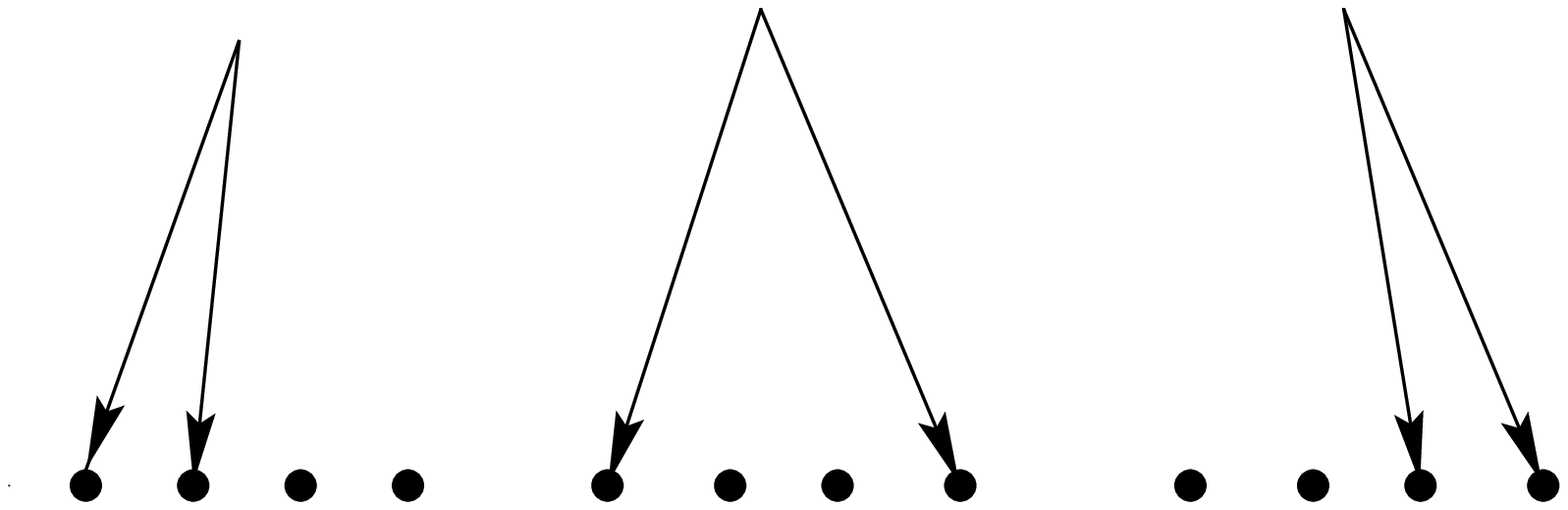}{10cm}

\vglue 1cm

\centerline{\bf Figure 1}

\centerline{The set $\G (2,1)$}

\medskip

We attach to any $\G \in \G (\ell , k)$ a polydifferential operator 
$\vp_{\G} (f_1 \ot \ldots \ot f_{\ell + 2k})$ as follows. Let $\xi_1 , 
\ldots 
\xi_{\ell}$ be vector fields. We set:
$$
\vp_{\G} (f_1 \ot \ldots \ot f_{\ell + 2k}) = \build{\rm Alt}_{\xi_1 \ldots 
\xi_{\ell}}^{} \prod_{i=1}^{\ell + 2k} \vp_{\G}^i (f_i) \eqno (26)
$$
where
$$
\vp_{\G}^i (f) = \left\{ \matrix{
f \, , &\hbox{if the $i$-th vertex of the} \hfill \cr
&\hbox{second type is not an endpoint} \hfill \cr
&\hbox{of an edge in the graph} \ \G \hfill \cr
\cr
\xi_j (f) \, , &\hbox{if the $i$-th vertex of the} \hfill \cr
&\hbox{second type is the $j$-th endpoint} \hfill \cr
&\hbox{(from left to right)} \hfill \cr
} \right.
$$

We set $(\g \in T_{\rm poly}^{\ell -1})$:
$$
\wt{\vp}_{\rm HKR}^{\rm cycl} (\g \ot u^k) = {k! \over (\ell + 2k)!} 
\sum_{\G 
\in \G (\ell , k)} \vp_{\G} \eqno (27)
$$
and
$$
\vp_{\rm HKR}^{\rm cycl} = {1 \over \ell + 2k + 1} \sum \wt{\vp}_{\rm 
HKR}^{\rm cycl} (f_1 \ot \ldots \ot f_{\ell + 2k}) \eqno (28)
$$
where $\sum = 1+C+\ldots + C^{\ell + 2k}$.

One can see that the definition (28) coincides with Examples 2.2 (i) -- (v) 
in the simplest cases.

\medskip

\noindent {\bf Theorem.} {\it The map $\vp_{\rm HKR}^{\rm cycl} : T_{\rm 
poly}^{\bu} \ot \Cb \, [u] \ra [\Dc_{\rm poly}^{\bu}]_{\rm cycl}$ is a map of 
the complexes
$$
\vp_{\rm HKR}^{\rm cycl} : \{ T_{\rm poly}^{\bu} \ot \Cb \, [u] , d_{\rm div} 
\} 
\ra \{ [\Dc_{\rm poly}^{\bu}]_{\rm cycl} , d_{\rm Hoch} \}
$$
and also a quasiisomorphism of the complexes.}

\medskip

We prove this Theorem in Sections 2.4 -- 2.6.

\medskip

\noindent {\bf 2.4.} To prove Theorem 2.3 we need some preparations.

\medskip

\noindent {\bf 2.4.1. Lemma.} {\it For any $k>0$, $\ell \geq 0$, and $\G 
\in 
\G (\ell , k)$, the cochain $\vp_{\G} (f_1 \ot \ldots \ot f_{\ell+2k})$ is a 
Hochschild coboundary.}

\medskip

\noindent {\it Proof.} For the graph $\G \in \G (\ell , k)$, $k > 0$, we 
define a graph $\wt \G$ as follows: the graph $\wt \G$ has 1 vertex of the 
first type, $\ell + 2k - 1$ vertices of the second type, and $\ell$ edges. 
We just short the first maximal sequence of consecutive free ($=$ not 
endpoints of edges) vertices of second type on 1 vertex, see Figure~2

\vglue 1cm

\sevafig{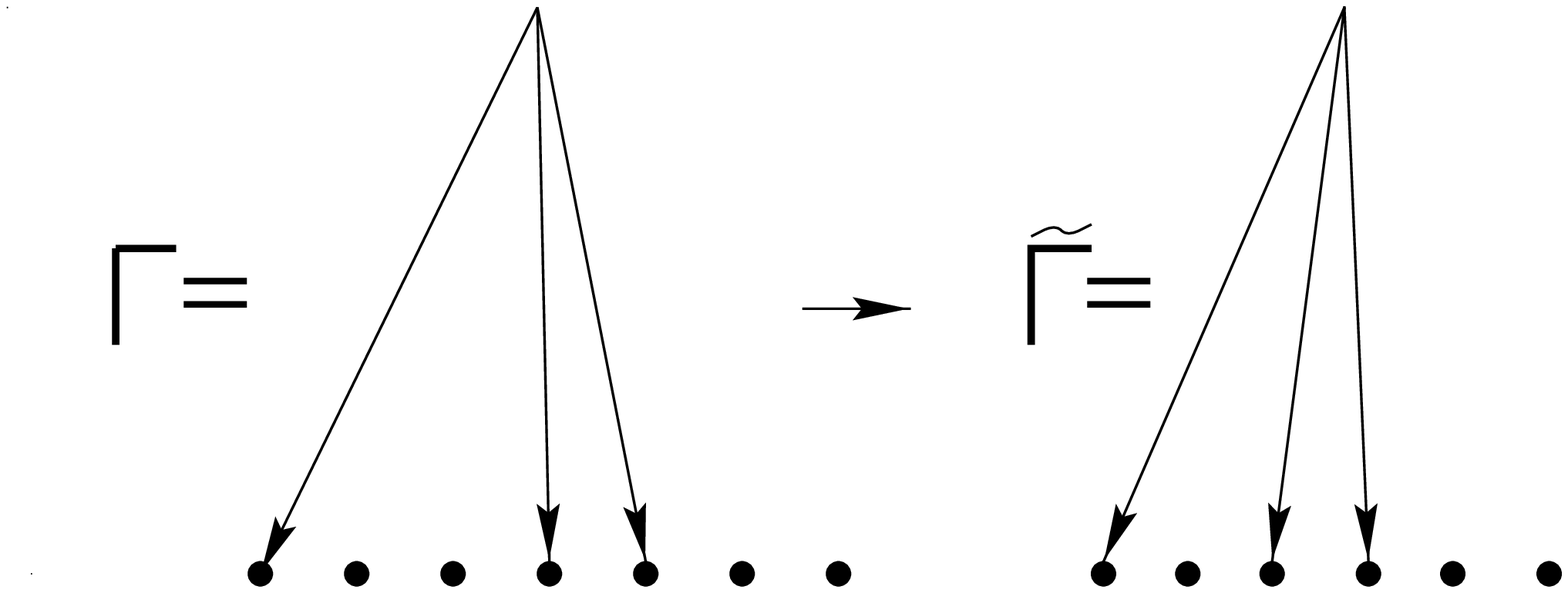}{10cm}

\vglue 1cm

\centerline{\bf Figure 2}

\centerline{$\G \in \G (3,2)$ and $\wt \G$}

\medskip

We define the cochain $\vp_{\wt{\G}} (f_1 \ot \ldots \ot f_{\ell+2k-1})$ 
analogously to (26). It is easy to see that
$$
(d_{\rm Hoch} \, \vp_{\wt{\G}}) (f_1 \ot \ldots \ot f_{\ell + 2k}) = 
\vp_{\G} 
(f_1 \ot \ldots \ot f_{\ell + 2k}) \, . 
$$
\hfill \xx

\bigskip

\noindent {\bf Corollary.} {\it For any $k,\ell \geq 0$, and any $\G \in \G 
(\ell , k)$ the cochain $\vp_{\G} (f_1 \ot \ldots \ot f_{\ell + 2k})$ is a 
Hochschild cocycle.} \hfill \xx

\bigskip

\noindent {\bf 2.4.2.} Here we study in which sense the map $\sum = 1+C+C^2 + 
\ldots + C^{i+1} : \Dc_{\rm poly}^i \ra \Dc_{\rm poly}^i$ is compatible with 
the Hochschild differential $d_{\rm Hoch}$. The results of this Subsection 
hold for any associative algebra $A$ with a 
trace functional $\int : A \ra \Cb$, i.e. $\int a \cdot b = \int b \cdot a$ 
for any $a,b \in A$, provided the condition that for any cochain $\psi (a_1 
\ot \ldots \ot a_k)$ there exists the unique cochain $(C\psi) (a_1 \ot 
\ldots \ot a_k)$ such that
$$
\int \psi (a_1 \ot \ldots \ot a_n) \cdot a_{n+1} = (-1)^n \int (C\psi) (a_2 
\ot \ldots \ot a_{n+1}) \cdot a_1 \, .
$$

\smallskip

\noindent {\bf Lemma.} {\it For $\psi \in {\rm Hom}_{\Cb} (A^{\ot k} , A)$ 
one has
$$
\eqalign{
& \ \int f_{k+2} \cdot (d_{\rm Hoch} \sum \psi - \sum d_{\rm Hoch} \, \psi) 
(f_1 \ot \ldots \ot f_{k+1}) = \cr
= & \ (-1)^{k-1} \int \phi (f_1 \ot \ldots \ot f_{n+2}) \cr
}
$$
where}
$$
\eqalignno{
& \ \phi (f_1 \ot \ldots \ot f_{n+2}) = \cr
= & \ \psi (f_1 \ot \ldots \ot f_n) \cdot f_{n+1} \cdot f_{n+2} + 
(-1)^{k-1} \, f_1 \cdot \psi (f_2 \ot \ldots \ot f_{n+1}) \cdot f_{n+2} + \cr
+ & \ f_1 \cdot f_2 \cdot \psi (f_3 \ot \ldots \ot f_{n+2}) + (-1)^{k-1} \, 
f_2 \cdot f_3 \cdot \psi (f_4 \ot \ldots \ot f_{n+2} \ot f_1) \cr
+ & \ \ldots \pm f_n \cdot f_{n+1} \cdot \psi (f_{n+2} \ot f_1 \ot \ldots 
\ot 
f_{n-1}) \, . &(29) \cr
}
$$

\smallskip

\noindent {\it Proof.} It is a straightforward calculation. \hfill \xx

\bigskip

\noindent {\bf 2.5.} Here we prove that the map
$$
\vp_{\rm HKR}^{\rm cycl} : T_{\rm poly}^{\bu} \ot \Cb \, [u] \ra [\Dc_{\rm 
poly}^{\bu}]_{\rm cycl}
$$
is a morphism of the complexes.

\smallskip

We have to prove that
$$
\vp_{\rm HKR}^{\rm cycl} ({\rm div} (\g) \ot u^{k+1}) = d_{\rm Hoch} \, 
\vp_{\rm HKR}^{\rm cycl} (\g \ot u^k) \eqno (30)
$$
for any $\g = \xi_1 \wdg \ldots \wdg \xi_{\ell} \in T_{\rm poly}^{\ell - 
1}$, 
and any $k$. We have:
$$
\eqalign{
& \ \int f_{\ell + 2k +2} \cdot d_{\rm Hoch} \, \vp_{\rm HKR}^{\rm cycl} 
(\g \ot u^k) (f_1 \ot \ldots \ot f_{\ell + 2k + 1}) \cdot \Om = \cr
= & \ \int f_{\ell + 2k +2} \cdot \left( d_{\rm Hoch} \cdot {1 \over \ell 
+2k 
+1} \sum \wt{\vp}_{\rm HKR}^{\rm cycl} \right) (f_1 \ot \ldots \ot f_{\ell + 2k 
+ 
1}) \Om = \cr
}
$$
by Lemma 2.4.2
$$
\eqalign{
& \ \int f_{\ell + 2k +2} \cdot \left( {1 \over \ell +2k +1} \sum d_{\rm 
Hoch} \right) (f_1 \ot \ldots \ot f_{\ell + 2k + 1}) \cdot \Om \cr
& +\int \phi(f_1\ot\ldots\ot f_{\ell+2k+2})\cdot\Om= \cr
}
$$
by  Corollary 2.4.1
$$
\int \phi (f_1 \ot \ldots \ot f_{\ell + 2k + 2}) \cdot \Om \, .
$$

On the other hand, for any $\G \in \G (\ell - 1 , k + 1)$ and any vector fields 
$\xi_1 , \ldots , \xi_{\ell}$ we define a cochain $\ov{\vp}_{\G} (f_1 \ot 
\ldots \ot f_{\ell + 2k + 2})$ as follows:
$$
\ov{\vp}_{\G} (f_1 \ot \ldots \ot f_{\ell + 2k + 2}) = \sum_{
\matrix{ 
\hbox{over all free} 
\cr 
\hbox{vertices $j$ of} 
\cr 
\hbox{second type and} 
\cr 
j = \ell + 2k + 2}} \ov{\vp}_{\G}^{(j)} (f_1 \ot \ldots \ot f_{\ell + 
2k 
+ 2}) \eqno (31)
$$
(free $=$ not endpoint) where
$$
\ov{\vp}_{\G}^{(j)} (f_1 \ot \ldots \ot f_{\ell + 2k + 2}) = \build{\rm 
Alt}_{\xi_1 , \ldots ,\xi_{\ell}}^{} \prod_{i=1}^{\ell + 2k +2} 
\vp_{\G}^{(j),i} (f_i) \eqno (32)
$$
and
$$
\vp_{\G}^{(j),i} (f) = \left\{ \matrix{
f \, , \hfill &\hbox{if} &i \ \hbox{is a free vertex of} \ \G \ \hbox{and} 
\ 
i \ne j \hfill \cr
\xi_{\ell} (f) \hfill &\hbox{if} &i=j \hfill \cr
\xi_s (f) \hfill &\hbox{if} &i \ \hbox{is $s$-th endpoint.} \hfill \cr
} \right.
$$
We set
$$
\ov{\vp} = \sum_{\G\in\G(\ell-1,k+1)} \ov{\vp}_{\G} \, . \eqno (33)
$$

\smallskip

\noindent {\bf 2.5.1. Lemma.} {\it Let} $\g = \xi_1 \wdg \ldots \wdg 
\xi_{\ell}$
$$
\eqalign{
& \ \int \ov{\vp}_{\ell , k} (f_1 \ot \ldots \ot f_{\ell + 2k + 2}) \cdot 
\Om = \cr
= & \ \pm \int \sum_{\G \in \G (\ell - 1 , k+1)} \vp_{\G} (({\rm div} \, 
\g) 
\ot u^{k+1}) (f_1 \ot \ldots \ot f_{\ell + 2k + 1}) \cdot f_{\ell + 2k + 2} 
\cdot \Om \, . \cr
}
$$

\smallskip

\noindent {\it Proof.} It is clear. \hfill \xx

\bigskip

\noindent {\bf 2.5.2. Lemma.} {\it Let $\phi (f_1 \ot \ldots \ot f_{\ell + 
2k 
+ 2})$ is defined from the cochain
$$
\sum_{\G \in \G (\ell , k)} \vp_{\G} (f_1 \ot \ldots \ot f_{\ell + 2k})
$$
by $(29)$. Then}
$$
\eqalign{
& \ (\ell +2k +1) \cdot \int \phi (f_1 \ot \ldots \ot f_{\ell + 2k +2}) 
\cdot \Om = \cr
= & \ \pm k \cdot \int \ov{\vp}_{\ell , k} (f_1 \ot \ldots \ot f_{\ell + 2k 
+2}) \cdot \Om \, . \cr
}
$$

\smallskip

\noindent {\it Proof.} It is straightforward. \hfill \xx

\bigskip

Now (30) follows from Lemma 2.5.1 and Lemma 2.5.2.

\medskip

\noindent {\bf 2.6.} We have proved that $\vp_{\rm HKR}^{\rm cycl} : T_{\rm 
poly}^{\bu} \ot \Cb \, [u] \ra [\Dc_{\rm poly}^{\bu}]_{\rm cycl}$ is a map of 
the complexes. It remains to prove that it is a quasiisomorphism. 

\smallskip

We know that both complexes have the same cohomology:
$$
\eqalign{
& \ H^i (T_{\rm poly}^{\bu} \ot \Cb \, [u] , d_{\rm div}) = H^i ([ \Dc_{\rm 
poly}^{\bu}]_{\rm cycl} , d_{\rm Hoch}) \cr \cr
= & \ \left\{ \matrix{
[T_{\rm poly}^i]_{\rm div} , &i \leq d-1 \hfill \cr
\Cb , \hfill &i = d+2k-1 , &k \geq 0 \, . \cr
} \right. \cr
}
$$
We know that $\vp_{\rm HKR}^{\rm cycl} \mid_{[T_{\rm poly}^{\bu}]_{\rm div}} = 
\vp_{\rm HKR} \mid_{[T_{\rm poly}^{\bu}]_{\rm div}}$, and therefore the map
$$
\vp_{\rm HKR}^{\rm cycl} \mid_{[T_{\rm poly}^{\bu}]_{\rm div}} : [T_{\rm 
poly}^{\bu}]_{\rm div} \ra \{ [\Dc_{\rm poly}^{\bu}]_{\rm cycl} , d_{\rm Hoch} \} 
\, ,
$$
is an embedding on the level of cohomology. Only what remains to prove is that the 
map
$$
\vp_{\rm HKR}^{\rm cycl} : \Cb = H^{d+2k-1} (T_{\rm poly}^{\bu} \ot \Cb \, [u] , 
d_{\rm div}) \ra H^{d+2k-1} ([\Dc_{\rm poly}^{\bu}]_{\rm cycl} , d_{\rm Hoch})
$$
is an embedding. We assume that $\Om = dx_1 \wdg \ldots \wdg dx_d$. It is clear 
that $\t_k = \left( {\part \over \part x_1} \wdg \ldots \wdg {\part \over \part 
x_d} \right) \ot u^k$ is a cocycle in $\{ T_{\rm poly}^{\bu} \ot \Cb \, [u] , d_{\rm 
div} \}$, and it is not cohomologous to zero. We have to prove that $\vp_{\rm 
HKR}^{\rm cycl} (\t_k)$ is not cohomologous to zero for any $k \geq 0$.

\medskip

\noindent {\bf 2.6.1. Long exact sequence, associated with bicomplex 2.1.3.}

\smallskip

Let us denote by $C^{\bu \bu}$ the bicomplex 2.1.3. Then there exists the 
following short exact sequence of the bicomplexes:
$$
0 \ra C^{\bu \bu} [-2] \ra C^{\bu \bu} \ra (\hbox{first 2 columns}) \ra 0 \, .
$$
The cohomology of the last term are equal to the Hochschild cohomology, because 
the complex $K^{\bu}$ is acyclic (Lemma 2.1.2). We obtain the following long 
exact sequence:
$$
\ldots \ra H^{i-1} \ra HC^{i-2} \build \ra_{}^{S} HC^i \ra H^i \ra HC^{i-1} 
\build \ra_{}^{S} HC^{i+1} \ra H^{i+1} \ra \ldots
$$
where $H^{\bu}$ stands for the Hochschild cohomology and $HC^{\bu}$ stands for 
the cyclic cohomology (i.e., the total cohomology of bicomplex 2.1.3). It is 
clear from the Hochschild-Kostant-Rosenberg theorem that the map $S : HC^{i-2} 
\ra HC^i$ is an {\it isomorphism} for $i \geq d+1$ ($A = \Cb \, [x_1 , \ldots , 
x_d]$).

\bigskip

\noindent {\bf 2.6.2. Lemma.} {\it $S \, [\vp_{\rm HKR}^{\rm cycl} (\t_k)] = 
\mu \cdot 
[\vp_{\rm HKR}^{\rm cycl} (\t_{k+1})]$ where $\t_k =$ \break $\left( {\part \over 
\part x_1} \wdg \ldots \wdg {\part \over \part x_d} \right) \ot u^k$, where $\mu 
\not= 0$.}

\bigskip

It follows from this Lemma that $\vp_{\rm HKR}^{\rm cycl} (\t_k)$ is not 
cohomologous to $0$ for any $k \geq 0$, and that the map $\vp_{\rm HKR}^{\rm 
cycl} : \{ T_{\rm poly}^{\bu} \ot \Cb \, [u] , d_{\rm div} \} \ra \{ [\Dc_{\rm 
poly}^{\bu}]_{\rm cycl} , d_{\rm Hoch} \}$ is a quasiisomorphism of the 
complexes.

\medskip

\noindent {\it Proof of Lemma.} We assume that $\Om = dx_1 \wdg \ldots \wdg 
dx_d$. It is clear that ${\rm div} \left( {\part \over 
\part x_1} \wdg \ldots \wdg {\part \over \part x_d} \right) = 0$. We have to 
prove that $S \, [{\vp}_{\rm HKR}^{\rm cycl} (\t_k)] = \mu \cdot [{\vp}_{\rm 
HKR}^{\rm cycl} (\t_{k+1})]$. We prove here the last formula for $k=0$, the 
general case is analogous.

\smallskip

\noindent We consider ${\vp}_{\rm HKR}^{\rm cycl} (\t_0)$ as an element of the 
third column of the bicomplex 2.1.3. We have: $(1-C) ({\vp}_{\rm HKR}^{\rm 
cycl} (\t_0)) = 0$ and $d_{\rm Hoch} ({\vp}_{\rm HKR}^{\rm cycl} (\t_0)) = 0$. 
The rows of the bicomplex 2.1.3 are exact, and it follows from the first equation 
that there exists $\a$ such that $\sum \a = {\vp}_{\rm HKR}^{\rm cycl} (\t_0)$ 
(the element $\a$ lies in the second column of the bicomplex 2.1.3). We set $\b = 
d_K (\a)$ where $d_K$ is the differential in the complex $K^{\bu}$. We have:

\smallskip

\noindent (1) $d_K \, \b = d_K^2 \, \a = 0$

\smallskip

\noindent (2) $\sum \b = \sum d_K \, \a = d_{\rm Hoch} \sum \a = d_{\rm Hoch} \,
{\vp}_{\rm HKR}^{\rm cycl} (\t_0) = 0$.

\smallskip

\noindent We need to know explicit formulas for $\a$ and $\b$.

\smallskip

\noindent We can choose $\a = {1 \over k+1} \, {\vp}_{\rm HKR}^{\rm cycl} 
(\t_0)$.

\smallskip

\noindent Then, up to a constant factor,
$$
\b = d_K \, \a = \left\{ f_1 \ot \ldots \ot f_{d+1} \ra \build{\rm Alt}_{{\part 
\over \part x_1} , \ldots , {\part \over \part x_d}}^{} {\part \over \part x_1} 
\, (f_1) \cdot \ldots \cdot {\part \over \part x_d} \, (f_d) \cdot f_{d+1} 
\right\} \, .
$$
Furthermore, there exists $\g$ such that $(1-C) \, \g = \b$. We set $\t = d_{\rm 
Hoch} \, \g$. It is clear that $\t$ is a cocycle in the bicomplex, cohomologous 
to the element ${\vp}_{\rm HKR}^{\rm cycl} (\t_0)$ in the third column. The 
element $\t$ lies in the first column. We have to express $\t$ explicitly.

\smallskip

First of all, let us find $\g$ such that $(1-C) \, \g = \b$.

\smallskip

\noindent Denote by $\b^{(k)}$, $1 \leq k \leq d+1$, the cochain
$$
\b^{(k)} = \build{\rm Alt}_{{\part \over \part x_1} , \ldots , {\part \over \part 
x_d}}^{} {\part \over \part x_k} \, (f_1) \cdot \ldots \cdot {\part \over \part 
x_{k-1}} \, (f_{k-1}) \cdot f_k \cdot {\part \over \part x_k} \, (f_{k+1}) \cdot 
\ldots \cdot {\part \over \part x_d} (f_{d+1}) \, .
$$

\smallskip

\noindent {\bf Lemma.}
$$
C (\b^{(k)}) = \left\{ \matrix{
- \, \b^{(k-1)} +(-1)^{k+d-1} \, \b^{(d+1)} , \ k>1 \cr
\b^{(d+1)} , \ k=1 \, . \hfill \cr
} \right.
$$

\smallskip

\noindent {\it Proof.} It is a direct calculation. \hfill \xx

\bigskip

It follows from this lemma that, up to a nonzero constant, one can choose
$$
\g = \b^{(1)} - \b^{(2)} + \b^{(3)} - \b^{(4)} + \ldots \pm \b^{(d+1)} \, .
$$
Now it is clear that $\t = d_{\rm Hoch} \, \g$ is equal, up to a  nonzero constant factor, 
to ${\vp}_{\rm HKR}^{\rm cycl} (\t_1)$ (see formulas (27), (28)).

\smallskip

The case $k>0$ is analogous. \hfill \xx

\bigskip

We have proved that the map
$$
\vp_{\rm HKR}^{\rm cycl} : \{ T_{\rm poly}^{\bu} \ot \Cb \, [u] , d_{\rm div} \} \ra 
\{ [\Dc_{\rm poly}^{\bu}]_{\rm cycl} , d_{\rm Hoch} \}
$$
is a quasiisomorphism of the complexes.

\vglue 1cm

\noindent {\bf 3. Cyclic Formality morphism}

\smallskip

In this section we construct explicitly maps
$$
\matrix{
\vp_{\rm HKR}^{\rm cycl} &= &\Cc_1 : T_{\rm poly}^{\bu} \ot \Cb \, [u] \ra 
[\Dc_{\rm poly}^{\bu}]_{\rm cycl} \hfill \cr \cr
&&\Cc_2 : \Lb^2 (T_{\rm poly}^{\bu} \ot \Cb \, [u]) \ra [\Dc_{\rm 
poly}^{\bu}]_{\rm cycl} [-1] \hfill \cr \cr 
&&\Cc_3 : \Lb^3 (T_{\rm poly}^{\bu} \ot \Cb \, [u]) \ra [\Dc_{\rm 
poly}^{\bu}]_{\rm cycl} [-2] \hfill \cr
}
$$
and our main Conjecture states that these maps are the Taylor components of an 
$L_{\ify}$-morphism $\Cc : \{ T_{\rm poly}^{\bu} \ot \Cb \, [u] , d_{\rm div} 
\} \ra \{ [\Dc_{\rm poly}^{\bu}]_{\rm cycl} , d_{\rm Hoch} \}$, i.e. for any 
$\eta_1 , \ldots , \eta_n \in T_{\rm poly}^{\bu} \ot \Cb \, [u]$ and any 
functions $f_1 , \ldots , f_m$ one has:
$$
\eqalignno{
\sum_{i=1}^n \pm & \ \Cc_n (\eta_1 \wdg \ldots \wdg d_{\rm div} \, \eta_i 
\wdg \ldots \wdg \eta_n) (f_1 \ot \ldots \ot f_m) + \cr
+ & \ f_1 \cdot \Cc_n (\eta_1 \wdg \ldots \wdg \eta_n) (f_2 \ot \ldots \ot 
f_m) \cr
\pm & \ \Cc_n (\eta_1 \wdg \ldots \wdg \eta_n) (f_1 \ot \ldots \ot f_{m-1}) 
\cdot f_m \cr
+ & \ \sum_{i=1}^{m-1} \pm \, \Cc_n (\eta_1 \wdg \ldots \wdg \eta_n) (f_1 \ot 
\ldots \ot f_i f_{i+1} \ot \ldots \ot f_m) + \cr
+ & \ \sum_{i \ne j} \pm \, \Cc_{n-1} ([\eta_i , \eta_j] \wdg \eta_1 \wdg 
\ldots \wdg \wh{\eta}_i \wdg \ldots \wdg \wh{\eta}_j \wdg \ldots \eta_n) 
(f_1 \ot \ldots \ot f_m) + \cr
+ & \ {1 \over 2} \sum_{k,\ell \geq 1 \atop k + \ell = n} {1 \over k! 
\ell!} \sum_{\s \in \Si_n} \cr
\pm & \ [\Cc_k (\eta_{\s_1} \wdg \ldots \wdg  \eta_{\s_k}) , \Cc_{\ell} 
(\eta_{\s_{k+1}} \wdg \ldots \wdg \eta_{\s_n})] (f_1 \ot \ldots \ot f_m) = 
0 \, . &(34)
}
$$
(We refer reader to [K], Section 4 for the general theory of 
$L_{\ify}$-algebras).

The idea of the construction of the maps $\Cc_1 , \Cc_2 , \Cc_3 , \ldots$ 
is the following.

We consider $[T_{\rm poly}^{\bu}]_{\rm div} \ot \Cb \, [u]$ as a $dg$ Lie 
subalgebra in $ \{ T_{\rm poly}^{\bu} \ot \Cb \, [u] , d_{\rm div} 
\} $ (with zero differential), 
and construct (with proofs) an $L_{\ify}$-morphism
$$
\wt{\Cc} : [T_{\rm poly}^{\bu}]_{\rm div} \ot \Cb \, [u] \ra \Dc_{\rm 
poly}^{\bu} \, .
$$
Here the idea is very close to [K].

In the next step, we apply the map $\sum$ to $\Dc_{\rm poly}^{\bu}$. The 
map $\sum$ is not a map of $dg$ Lie algebras (see Lemma 2.4.2), but it 
turns out that the composition
$$
\left[ \sum \right] \circ \wt{\Cc} : [T_{\rm poly}^{\bu}]_{\rm div} \ot \Cb \, 
[u] \ra [\Dc_{\rm poly}^{\bu}]_{\rm cycl}
$$
where $\left[ \sum \right] : \Dc_{\rm poly}^i \ra [\Dc_{\rm 
poly}^i]_{\rm cycl}$ is equal to
$$
\left[ \sum \right] = {1 \over i+2} \, \sum
$$
is still an $L_{\ify}$-morphism (it is a conjecture), and, moreover, the 
same formulas define an $L_{\ify}$-morphism $ \{ T_{\rm poly}^{\bu} \ot \Cb \, [u] , d_{\rm div} 
\} \ra \{ [\Dc_{\rm poly}^{\bu}]_{\rm cycl} , d_{\rm Hoch} \}
$. We have checked this fact only for 
$\Cc_1 = \vp_{\rm HKR}^{\rm cycl}$ in Section 2, and this check is quite 
nontrivial.

\medskip

\noindent {\bf 3.1. $L_{\ify}$-map $\wt{\Cc} : [T_{\rm poly}^{\bu}]_{\rm 
div} \ot \Cb \, [u] \ra \Dc_{\rm poly}^{\bu}$}

\medskip

\noindent {\bf 3.1.1. Admissible graphs}

\medskip

\noindent {\bf Definition.} Admissible graph $\G$ is an oriented graph with 
labels and two types of edges: usual edges, and {\it dashed edges}, such that

\smallskip

\item{1)} the set of vertices $V_{\G}$ is $\{ 1, \ldots , n \} \sqcup \{ 
\ov 1 , \ldots , \ov m \}$ where $n,m \in \Zb_{\geq 0}$; vertices from the 
set $\{ 1, \ldots , n \}$ are called vertices of the first type, vertices 
from $\{ \ov 1 , \ldots , \ov m \}$ are called vertices of the second type,

\smallskip

\item{2)} every edge $(v_1 , v_2 ) \in E_{\G}$ starts at a vertex of the 
first type, $v_1 \in \{ 1, \ldots , n \}$,

\smallskip

\item{3)} there are no loops,

\smallskip

\item{4)} dashed edges end only on the vertices of second type and appear 
only in pairs, i.e. it is a pair of edges $(v , \ov i)$ and $(v,\ov{i+1})$ 
where $v \in \{ 1, \ldots , n \}$ and $\ov i$, $\ov{i+1}$ are elements from 
$\{ \ov 1 , \ov 2 , \ldots , \ov m \}$; we will speak about a {\it dashed pair},

\smallskip

\item{5)} for every vertex $\ell \in \{ 1, \ldots , n \}$ the set of usual 
(not dashed) edges
$$
{\rm Star} \, (\ell) := \{ (v_1 , v_2) \in E_{\G} \mid v_1 = \ell \}
$$
is labelled by symbols $(e_k^1 , \ldots , e_k^{\# \, {\rm Star} (\ell)})$.

\smallskip

A typical admissible graph is shown on Figure~3

\vglue 1cm

\sevafig{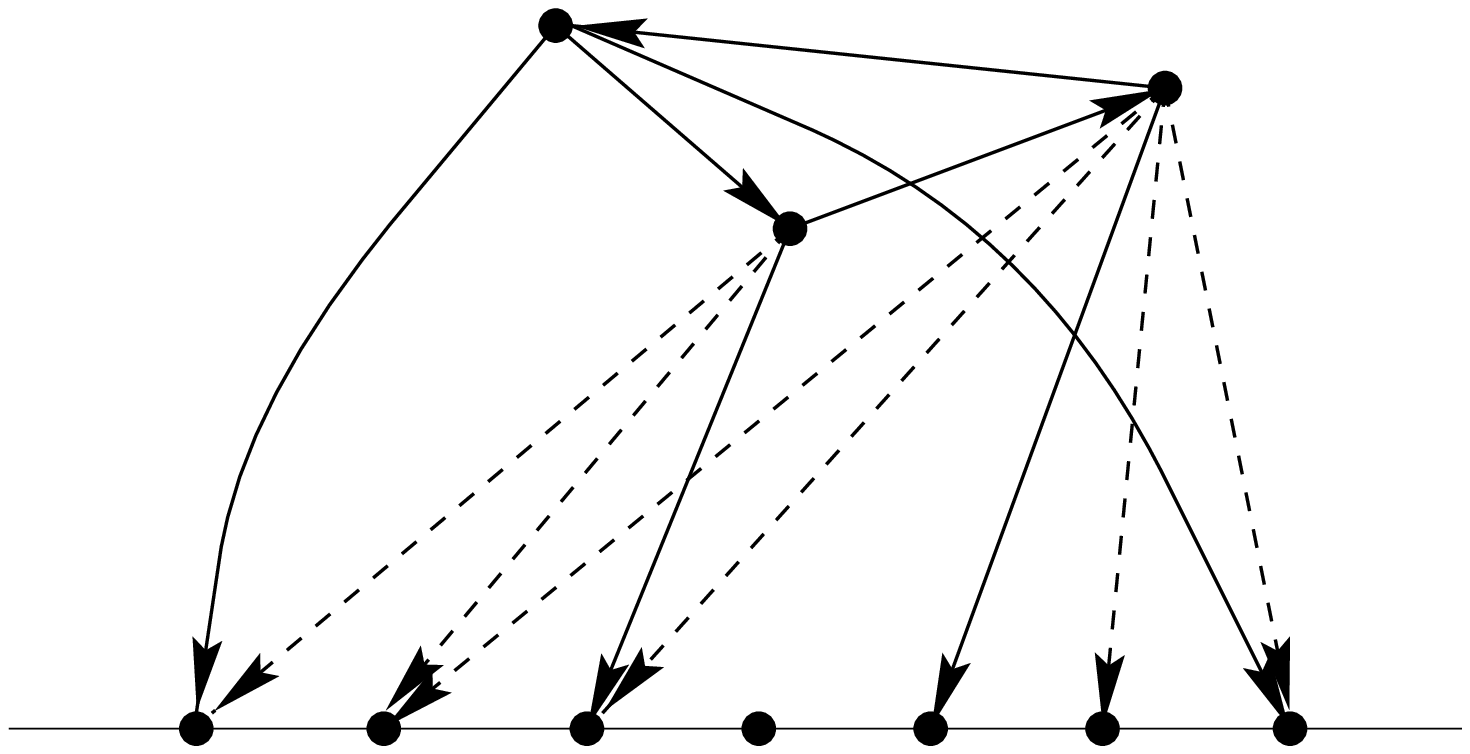}{10cm}

\vglue 1cm

\centerline{\bf Figure 3}

\centerline{An admissible graph with 3 dashed pairs}

\bigskip

\noindent {\bf 3.1.2. Configuration spaces}

\smallskip

We work with the same configuration spaces as in [K]. Let us recall the 
definitions.

Let $n,m$ be non-negative integers satisfying the inequality $2n+m \geq 2$. 
We denote by Conf$_{n,m}$ the product of configuration space of the upper 
half-plane with the configuration space of real line:
$$
\eqalign{
{\rm Conf}_{n,m} = & \ \{ (p_1 , \ldots , p_n ; q_1 , \ldots , q_m) \mid 
p_i \in \Hc \, , \ q_j \in \Rb \, , \cr
& \ p_{i_1} \ne p_{i_2} \ \hbox{for} \ i_1 \ne i_2 \ \hbox{and} \ q_{j_1} 
\ne q_{j_2} \ \hbox{for} \ j_1 \ne j_2 \} \, . \cr
}
$$
The group
$$
G^{(1)} = \{ z \mpo az + b \mid a,b \in \Rb \, , \ a > 0 \}
$$
acts on the space Conf$_{n,m}$. It follows from the condition $2n+m \geq 2$ 
that this action is free. The quotient space $C_{n,m} = {\rm Conf}_{n,m} / 
G^{(1)}$ is a smooth manifold of dimension $2n+m-2$.

Analogously, we introduce simpler spaces Conf$_n$ and $C_n$ for any $n \geq 
2$:
$$
{\rm Conf}_n = \{ (p_1 , \ldots , p_n) \mid p_i \in \Cb \, , \ p_i \ne p_j 
\ \hbox{for} \ i \ne j \}
$$
$$
C_n = {\rm Conf}_n  / G^{(2)} \, , \ {\rm dim} \, C_n = 2n-3
$$
where $G^{(2)}$ is a group
$$
G^{(2)} = \{ z \mpo az + b \mid a \in \Rb \, , \ b \in \Cb \, , \ a > 0 \} \, 
.
$$

We construct (following [K], Section 5) compactifications ${\ov C}_{n,m}$ 
of $C_{n,m}$ (and compactifications ${\ov C}_n$ of $C_n$). These 
compactifications are ``manifolds with corners''.

Let us describe the strata of codimension 1. There are two types of strata 
of codimension 1 in ${\ov C}_{A,B}$:

\smallskip

\item{S1)} points $p_i \in \Hc$ for $i \in S \sbsq A$ where $\# \, S \geq 
2$ move close to each other and far from $\Rb$, the corresponding boundary 
stratum
$$
\part_S \, \ov{C}_{A,B} = C_S \ts C_{(A \bsh S)} \sqcup \{ pt \} \, , \ B
$$

\item{S2)} points $p_i \in \Hc$ for $i \in S \sbsq A$ and points $q_j \in 
\Rb$ for $j \in S' \sbsq B$, where $2 \, \# \, S + \# \, S' \geq 2$,  all 
move close to each other and to $\Rb$, with at least one point left outside 
$S$ and $S'$, i.e. $\# \, S + \# \, S' \leq \# \, A + \# \, B-1$. The 
corresponding boundary stratum is
$$
\part_{S,S'} \, \ov{C}_{A,B} = C_{S,S'} \ts C_{A \bsh S} \, , \ (B \bsh S') 
\sqcup \{ pt \} \, .
$$

It is instructional to draw low-dimensional spaces $C_{n,m}$. The simplest 
one, $C_{1,0} = {\ov C}_{1,0}$ is just a point. The space $C_{0,2} = {\ov 
C}_{0,2}$ is a two-element set. The space $C_{1,1}$ is an open interval, 
and its closure ${\ov C}_{1,1}$ is a closed interval.

The space $C_{2,0}$ is diffeomorphic to $\Hc \bsh \{ 0 + 1 \cdot i \}$. The 
closure ${\ov C}_{2,0}$ is shown on Fig.~4.

\vglue 1cm

\sevafig{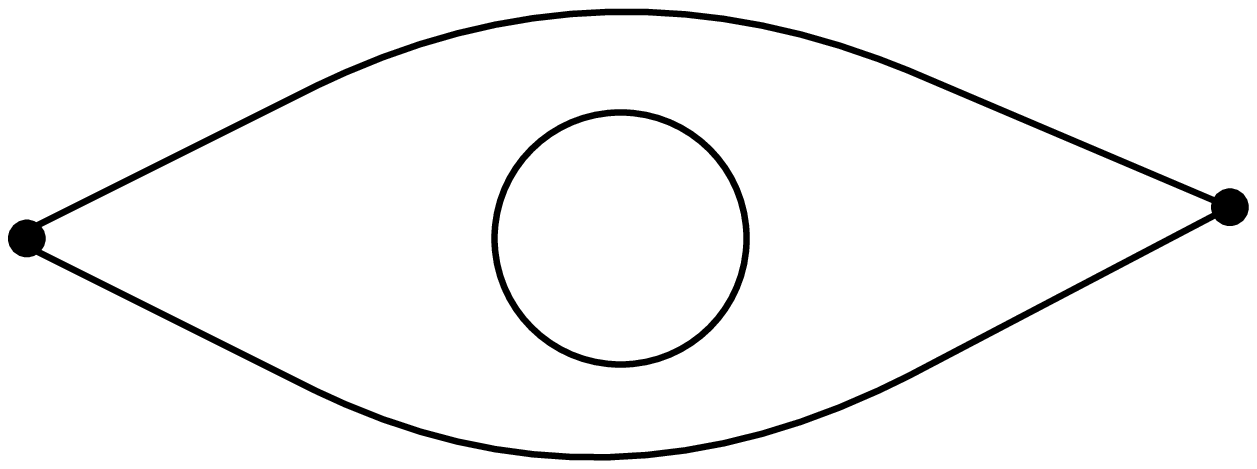}{5cm}

\vglue 1cm

\centerline{\bf Figure 4}

\centerline{The space ${\ov C}_{2,0}$ (``the Eye'')}

\medskip

\noindent See [K], Section 5 for more details.

\medskip

\noindent {\bf 3.1.3. Differential forms on configuration spaces}

\smallskip

The space ${\ov C}_{2,0}$ is homotopy equivalent to the standard circle 
$S^1 \sm \Rb / 2\pi \, \Zb$. Moreover, one of its boundary components, the 
space $C_2 = {\ov C}_2$ is naturally $S^1$. The other component of the 
boundary is the union of two closed intervals (copies of ${\ov C}_{1,1}$) 
with identified endpoints.

\medskip

\noindent {\bf Definition.} An angle map is a smooth map $\phi : {\ov 
C}_{2,0} \ra \Rb / 2\pi \, \Zb$ such that the restriction of $\phi$ to $C_2 
\sm S^1$ is the angle measure in anti-clockwise direction from the vertical 
line, and $\phi$ maps the whole upper interval ${\ov C}_{1,1} \sm [0,1]$ of 
${\ov C}_{2,0}$ to a point in $S^1$.

\medskip

We denote by $G_{n,m,2k}$ the set of all the admissible graphs with $n$ 
vertices of the first type, $m$ vertices of the second type, $k$ dashed 
pairs and $2n+m-2k-2$ usual edges. Let $\G \in G_{n,m,2k}$.

We define the weight $W_{\G}$ of the graph $\G$ by the following formula:
$$
W_{\G} = \prod_{\ell = 1}^n (k_{\ell})! \cdot \prod_{\ell = 1}^n {1 \over 
(\# \, {\rm Star} (\ell))!} \cdot {1 \over (2\pi)^{2n+m-2}} \int_{{\ov 
C}_{n,m}^+} \bigwedge_{e \in E_{\G}} d \, \vp_e \, . \eqno (35)
$$

Let us explain what is written here. The domain of integration ${\ov 
C}_{n,m}^+$ is a connected component of ${\ov C}_{n,m}$ which is the 
closure of configurations for which points $q_j$, $1 \leq j \leq m$ on 
$\Rb$ are placed in the increasing order: $q_1 < \ldots < q_m$. Every edge 
$e$ of $\G$ defines a map from ${\ov C}_{n,m}$ to ${\ov C}_{2,0}$ or to 
${\ov C}_{1,1} \sbs {\ov C}_{2,0}$ (we consider ${\ov C}_{1,1} \sbs {\ov 
C}_{2,0}$ as the lower interval of the Eye). The pull-back of the angle 
function $\phi$ by the map ${\ov C}_{n,m} \ra {\ov C}_{2,0}$ corresponding 
to edge $e$ is denoted by $\phi_e$. The number $k_{\ell}$ is the number of 
dashed pairs starting at the vertice $\ell$. Finally, the ordering in the 
wedge product of 1-forms $d\phi_e$ is fixed by enumeration of the set of 
sources of edges and by the enumeration of the set of edges with a given 
source (the dashed pairs can be counted in any order).

\medskip

\noindent {\bf 3.1.4. Pre-$L_{\ify}$-morphism associated with graphs}

\smallskip

For any admissible graph $\G \in G_{n,m,2k}$ (it has $n$ vertices of the 
first type, $m$ vertices of the second type, $k$ dashed pairs, and 
$2n-2k+ m - 2$ usual edges) we define a linear map $\wt{\Cc}_{\G} : 
\bigotimes^n ([T_{\rm poly}^{\bu}]_{\rm div} \ot \Cb \, [u]) \ra \Dc_{\rm 
poly}^{\bu} [1-n]$. This map has only one nonzero component 
$(\wt{\Cc}_{\G})_{(\ell_1 , k_1), \ldots , (\ell_n , k_n)}$ ($(\ell_j , 
k_j)$ here stands for $[T_{\rm poly}^{\ell_j}]_{\rm div} \ot u^{k_j}$),
 where $\ell_j 
= \# \, {\rm Star} \, (j) - 1$, and $k_j$ is the number of dashed pairs 
starting at the vertex $j$ of the first type.

Let $\eta_1 , \ldots , \eta_n \in [T_{\rm poly}^{\bu}]_{\rm div} \ot \Cb \, 
[u]$, $\eta_j = \g_j \ot u^{k_j}$ and $\g \in [T_{\rm poly}^{\ell_j}]_{\rm 
div}$, and let $f_1 , \ldots , f_m$ be functions on $\Rb^d$. We are going 
to write a formula for function $\Phi$ on $\Rb^d$:
$$
\Phi = \wt{\Cc}_{\G} (\eta_1 \ot \ldots \ot \eta_n) (f_1 \ot \ldots \ot 
f_m) \, .
$$
The formula for $\Phi$ is the sum over all configurations of indices 
running from $1$ to $d$, labeled by $ \ov{E}_{\G}$ where  $ \ov{E}_{\G}$ 
is the set of usual (not dashed) edges of the 
graph $\G$:
$$
\Phi = \sum_{I : \ov{E}_{\G} \{ 1 , \ldots , d \}} \Phi_I \eqno (36)
$$
where $\Phi_I$ is the product over all $n+m$ vertices of $\G$ of certain 
partial derivatives of functions $f_j$ and of polyvector fields $\g_i$.

Namely, with each vertex $i$, $1 \leq i \leq n$ of the first type we 
associate function $\psi_i$ on $\Rb^d$, where
$$
\psi_i = \lgl \g_i , dx^{I (e_i^1)} \ot \ldots \ot dx^{I (e_i^{\ell_i +1})} 
\rgl \, .
$$
Here we use the identification of polyvector fields with skew-symmetric 
tensor fields as
$$
\xi_1 \wdg \ldots \wdg \xi_{\ell + 1} \mpo \sum_{\s \in \sum_{\ell + 1}} 
{\rm sgn} (\s) \, \xi_{\s_1} \ot \ldots \ot \xi_{\s_{\ell +1}} \, .
$$
For each vertex $\ov j$ of the second type the associated function is 
defined $\psi_{\ov j}$ as $f_j$.

Now, at each vertex of graph $\G$ we put a function on $\Rb^d$ (i.e. 
$\psi_i$ of $\psi_{\ov j}$). Also, on edges of graph $\G$ there are indices 
$I(e)$ which label coordinates in $\Rb^d$. In the next step we put into 
each vertex $v$ instead of the function $\psi_v$ its partial derivative
$$
\left( \prod_{e \in \ov{E}_{\G} \, , \atop e = (*,v)} \part_{I(e)} \right) 
\psi_v \, , 
$$
and take the product over all vertices $v$ of $\G$. The result is by 
definition the summand $\phi_I$.

\medskip

\noindent {\it Remark.} The graphs we have considered in Section 2 are 
exactly graphs $\G \in G_{1,m,2\ell}$. Indeed, the dashed pair is by the 
definition a pair of edges $(1,\ov j )$ and $(1,\ov{j+1})$ for some $\ov j 
\in \{ \ov 1 , \ldots , \ov m \}$, and the graph $\G$ has $2\ell$ dashed 
edges and $m-2\ell$ usual edges. Every vertex of the second type is the 
endpoint of exactly one edge (usual or dashed), because otherwise the 
corresponding weight $W_{\G} = 0$.

\medskip

\noindent {\bf 3.1.5. $L_{\ify}$-morphism $\wt{\Cc} : [T_{\rm 
poly}^{\bu}]_{\rm div} \ot \Cb \, [u] \ra \Dc_{\rm poly}^{\bu}$}

\medskip

\noindent {\bf Theorem.} {\it The maps
$$
\wt{\Cc}_n = \sum_{m \geq 0} \sum_{k \geq 0} \sum_{\G \in G_{n,m,2k}} 
W_{\G} \cdot \wt{\Cc}_{\G} \eqno (37)
$$
are the Taylor components of an $L_{\ify}$-morphism
$$
\wt{\Cc} : [T_{\rm poly}^{\bu}]_{\rm div} \ot \Cb \, [u] \ra \Dc_{\rm 
poly}^{\bu} \, ,
$$
(or $\wt{\Cc} : \{[T_{\rm poly}^{\bu}] \ot \Cb \, [u],0\} \ra \Dc_{\rm 
poly}^{\bu}$).}

\smallskip

\noindent {\it Proof.} The proof is analogous to proof of the 
$L_{\ify}$-Formality conjecture in [K], Section 6. The left-hand side of 
(34) can be written as a linear combination
$$
\sum_{\G} C_{\G} \, \wt{\Cc}_{\G} (\eta_1 \ot \ldots \ot \eta_n) (f_1 \ot 
\ldots \ot f_m) \eqno (38)
$$
over admissible graphs with $n$ vertices of the first type, $m$ vertices of 
the second type, and $2n+m-3$ edges (usual and dashed).

Coefficients $C_{\G}$ in (38) are equal to quadratic-linear combinations of 
the weights $W_{\G'}$, $\G' \in G_{n,m,*}$.

We want to check that $C_{\G}$ vanishes for each $\G$.

The idea is to identify $C_{\G}$ with the integral over the boundary $\part 
\, \ov{C}_{n,m}$:
$$
\int_{\part \, \ov{C}_{n,m}} \bigwedge_{e \in E_{\G}} d \phi_e = 
\int_{\ov{C}_{n,m}} d \left( \bigwedge_{e \in E_{\G}} d \phi_e \right) = 0 
\, . \eqno (39)
$$
We have:
$$
0 = \int_{\part \, \ov{C}_{n,m}} \bigwedge_{e \in E_{\G}} d \phi_e = \sum_S 
\int_{\part_S \, \ov{C}_{n,m}} \bigwedge_{e \in E_{\G}} d \phi_e + 
\sum_{S,S'} \int_{\part_{S,S'} \, \ov{C}_{n,m}} \bigwedge_{e \in E_{\G}} d 
\phi_e \, . \eqno (40)
$$

The first summand in the r.h.s. of (40) does not vanish only for $\# \, S = 
2$ (see [K], Section 6.4.1), and this case corresponds to the summands in 
(34) with Schouten-Nijenhuis bracket. The second summand in the r.h.s. of 
(40), corresponds to summands in (34) with the Hochschild coboundary and 
with the Gerstenhaber bracket (see [K], Section 6.4.2). The only new thing 
is the dashed pairs.

\medskip

\noindent {\bf 3.1.5.1.}

\vglue 1cm

\sevafig{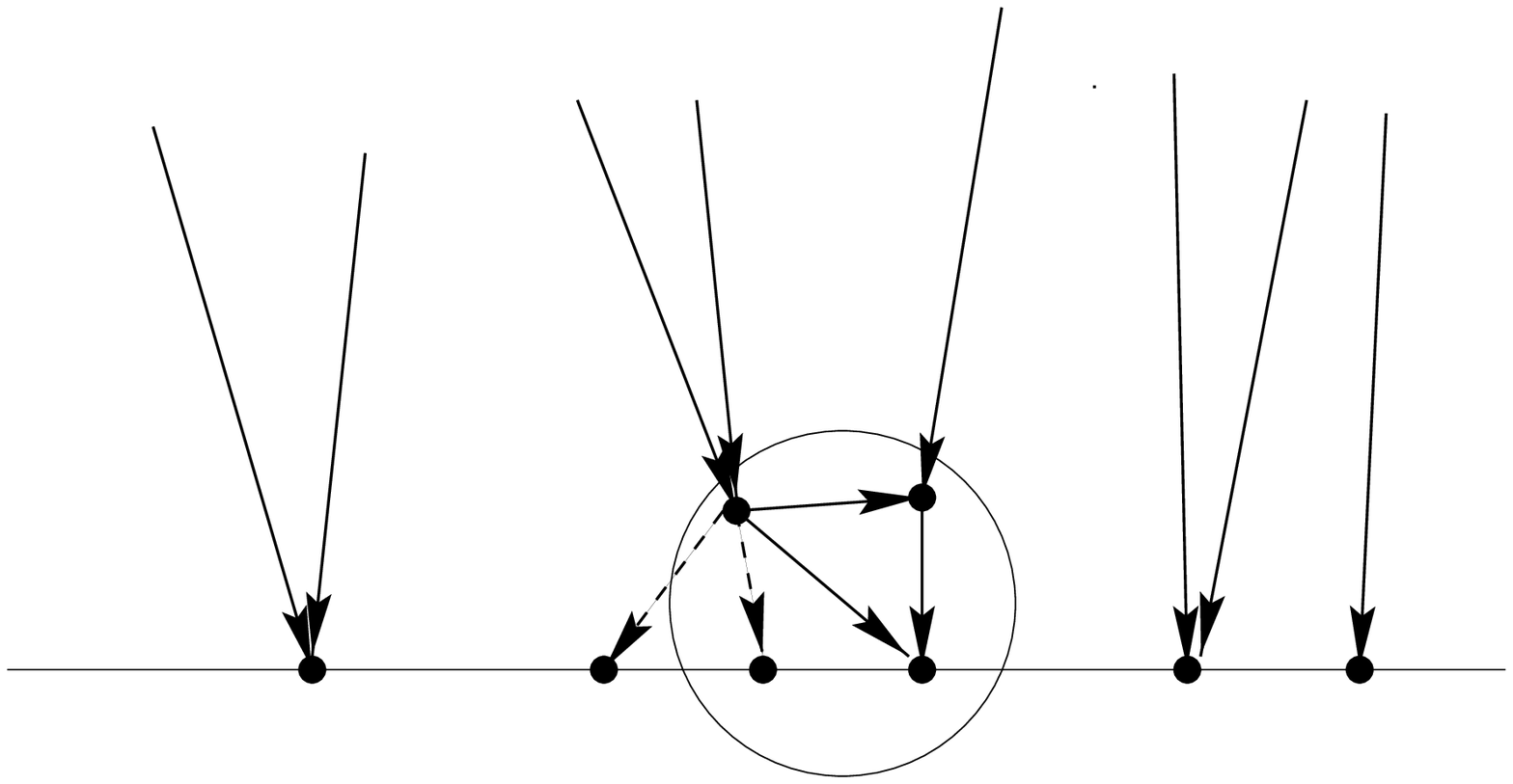}{10cm}

\vglue 1cm

\centerline{\bf Figure 5}

\medskip

When a vertex of the first type $v \in S$, in the boundary component 
$\part_{S,S'} \, \ov{C}_{n,m}$ the both endpoints of all the dashed pairs, 
starting at the vertex $v$, lie in $S'$; otherwise, the corresponding 
integral vanishes. The situation, shown on Figure 5, 
corresponds to zero integral.

\medskip

\noindent {\bf 3.1.5.2.} When two points of the first type with $k_1$ and 
$k_2$ dashed pairs move close to each other (case (S1) with $\# \, S = 2$) 
we obtain a vertex with $k_1 + k_2$ dashed pairs. The same final graph is 
corresponded to ${(k_1 + k_2)! \over k_1 ! \, k_2 !}$ graphs (i.e. we can 
select any $k_1$ dashed pairs from $k_1 + k_2$ dashed pairs as dashed pairs 
of the first vertex from two vertices which move close to each other). This 
is the cause of the appearance of the product $\build\prod_{\ell = 1}^{n} 
(k_{\ell})!$ in the formula (35) for weight $W_{\G}$. \hfill \xx

\bigskip

\noindent {\bf 3.2. The cyclic $L_{\ify}$-morphism}

\smallskip

We denote by $\left[ \sum \right]$ the operator $\build 
{1\over i+2}\cdot\sum_{}^{} 
:
 \Dc_{\rm poly}^i \ra [\Dc_{\rm poly}^i]_{\rm cycl}$, $i \geq -1$.

\medskip

\noindent {\bf Conjecture 1.} {\it The composition $\Cc = \left[ \sum 
\right] \circ \wt{\Cc}$ defines an $L_{\ify}$-morphism}
$$
\Cc : \{ [T_{\rm poly}^{\bu}]_{\rm div} \ot \Cb \, [u] , 0 \} \ra \{ [\Dc_{\rm 
poly}^{\bu}]_{\rm cycl} , d_{\rm Hoch} \} \, .
$$

\smallskip

\noindent {\bf Conjecture 2.} {\it The composition $\Cc = \left[ \sum 
\right] \circ \wt{\Cc}$ defines an $L_{\ify}$-morphism}
$$
\Cc : \{ T_{\rm poly}^{\bu} \ot \Cb \, [u] , d_{\rm div} \} \ra \{ [\Dc_{\rm 
poly}^{\bu}]_{\rm cycl} , d_{\rm Hoch} \} \, .
$$

\smallskip

Conjecture 2 is the cyclic Formality conjecture. We have proved (34) for 
$n=1$ and $\Cc = \left[ \sum \right] \circ \wt{\Cc}$ in Section 2. Let us 
note that it follows from Conjecture 2 and Theorem 2.3 that the 
$L_{\ify}$-morphism $\Cc$ is in fact an $L_{\ify}$-quasiisomorphism.

\vglue 1cm

\noindent {\bf 4. Globalization}

\smallskip

\noindent {\bf Corollary.} {\it Assuming Conjecture 2 the following is true. 
Let $M$ be a smooth manifold, $\Om$ be a 
volume form on $M$. Then there exists an $L_{\ify}$-quasiisomorphism}
$$
\Cc_M : \{ T_{\rm poly}^{\bu} (M) \ot \Cb \, [u] , d_{\rm div} \} \ra \{ 
[\Dc_{\rm poly}^{\bu} (M)]_{\rm cycl} , d_{\rm Hoch} \} \, .
$$

\smallskip

 The proof is analogous to the proof of the globalization of the 
$L_{\ify}$-Formality morphism ([K], Section 7), and we omit it here. 

\vglue 1cm

\noindent {\bf 5. Applications}

\smallskip

Here we consider first applications of the previous results to the 
deformation quantization.

\medskip

\noindent {\bf 5.1. Maurer-Cartan equation}

\smallskip

Let $g^{\bu}$ be a $dg$ Lie algebra.

\medskip

\noindent {\bf Definition.} An element $\g \in g^1$ satisfies the 
Maurer-Cartan equation iff $d \g + {1 \over 2} \, [\g , \g] = 0$.

\medskip

\noindent {\bf Lemma.} {\it Let 
$g_1^{\bu}$, $g_2^{\bu}$ be two $dg$ Lie algebras, and let $F : g_1^{\bu} 
\ra g_2^{\bu}$ be an $L_{\ify}$-morphism.
Let $\g \in g_1^1$ satisfies the Maurer-Cartan equation in $g_1^{\bu}$. Then
$$
F_1 (\g) + {1 \over 2} \, F_2 (\g , \g) + {1 \over 6} \, F_3 (\g , \g , \g) 
+ \ldots + {1 \over n!} \, F_n (\g , \ldots , \g) + \ldots = F(\g) \eqno 
(41)
$$
also satisfies the Maurer-Cartan equation in $g_2^{\bu} : d(F(\g)) + 
{1 \over 2} \, [F(\g) , F(\g)] = 0$.}

\medskip

\noindent {\it Example.} Let $\Uc : T_{\rm poly}^{\bu} \ra \Dc_{\rm 
poly}^{\bu}$ be an $L_{\ify}$-quasiisomorphism of M.~Kontsevich [K]. A 
solution of the Maurer-Cartan equation in $T_{\rm poly}^{\bu}$ is a 
bivector field $\g$ such that $[\g , \g] = 0$. Formula (41) produces a 
solution of the Maurer-Cartan equation in $\Dc_{\rm poly}^{\bu}$.

\medskip

\noindent {\bf Lemma.} {\it Let $A = C^{\ify} (\Rb^d)$. An element $\vp \in 
\Dc_{\rm poly}^1 (\Rb^d) \sbs {\rm Hom}_{\Cb} (A^{\ot 2} , A)$ satisfies 
the Maurer-Cartan equation iff $f*g = f \cdot g + \vp (f \ot g)$
 defines an associative star-product.} \hfill 
\xx

\bigskip

This Lemma and formula (41) produce a deformation quantization of the Poisson 
structure on 
$\Rb^d$, given by the bivector field $\g$.

\medskip

\noindent {\bf 5.2. Cyclically-invariant deformation quantization}

\smallskip

A solution of the Maurer-Cartan equation in the $dg$ Lie algebra $\{ T_{\rm 
poly}^{\bu} \ot \Cb \, [u] , d_{\rm div} \}$ is a bivector field $\g$ such 
that:
$$
\matrix{
(1) &[\g , \g] = 0 \hfill \cr \cr
(2) &{\rm div} (\g) = 0 \hfill \cr
} \eqno \matrix{(42) \cr \cr (43)}
$$
The cyclic $L_{\ify}$-Formality morphism
$$
\Cc : \{ T_{\rm poly}^{\bu} \ot \Cb \, [u] , d_{\rm div} \} \ra \{ [\Dc_{\rm 
poly}^{\bu}]_{\rm cycl} , d_{\rm Hoch} \}
$$
given in Conjecture 2 of Section 3.2 produces by formula (41) a solution 
$\psi$ of the Maurer-Cartan equation in $[\Dc_{\rm poly}^{\bu}]_{\rm 
cycl}$, i.e. an element $\psi : A^{\ot 2} \ra A$ such that
$$
\matrix{
(1) &f*g = f \cdot g + \psi (f \ot g) \ \hbox{is an} \hfill \cr
&\hbox{associative star-product on} \ \Rb^d \hfill \cr
} \eqno (44)
$$
$$
\matrix{
(2) &\int_{\Rb^d} (f*g) \cdot h \cdot \Om \ \hbox{is invariant with 
respect} \hfill \cr
&\hbox{to the cyclic permutation of} \ (f,g,h) \hfill \cr
} \eqno (45)
$$
for any functions $f,g,h \in C^{\ify} (\Rb^d)$ with compact support.

Moreover, the fact that $\wt{\Cc}$ is an $L_{\ify}$-quasiisomorphism 
(Section 2) coupled with standard deformation theory [K] allows us to deduce
 the following statement:

\medskip

\noindent {\bf Corollary.} {\it Assuming Conjecture 2 the following is true. 
Let $\Ac$ be the set of all Poisson 
structures on $\Rb^d$ satisfying $(43)$ modulo diffeomorphisms of $\Rb^d$ 
generating by vector fields $\xi$ such that $[ \g , \xi ] = 0$ and ${\rm 
div} \xi = 0$. Let $\Bc$ be the set of all star-products on $\Rb^d$
$$
f*g = f \cdot g + \hbar B_1 (f,g) + \hbar^2 B_2 (f,g) + \ldots
$$
which satisfy $(45)$ 
for any three functions $f,g,h$ with compact support modulo transformations
$$
f*g \ra T (T^{-1} (f) * T^{-1} (g))
$$
where $T(f) = f + \hbar T_1 (f) + \hbar^2 T_2 (f) + \ldots$ and 
$\int_{\Rb^d} T(f) \cdot g \cdot \Om = \int_{\Rb^d} f \cdot T^{-1} (g) \cdot \Om 
$ 
for 
any two functions $f,g$ on $\Rb^d$ with compact support.

Then the sets $\Ac$ and $\Bc$ are canonically isomorphic, and the 
isomorphism is given by formula $(41)$.} 

\bigskip

\noindent {\bf 5.3.} After the globalization of the cyclic 
$L_{\ify}$-Formality conjecture (Section 4) one can deduce from Conjecture 2
 the statement
analogous to Corollary 5.2 in the case of arbitrary $C^{\ify}$-manifolds.

\medskip

\noindent {\it Example.} Let $M$ be a symplectic manifold of dimension $2d$ 
with the symplectic form $\om$. Then there exists a deformation 
quantization of the algebra of smooth functions on this manifold such that 
for any 3 functions $f,g,h$ with compact support one has:
$$
\int_{\lgl M \rgl} (f*g) \cdot h \cdot \om^d = \int_{\lgl M \rgl} (g*h) 
\cdot f \cdot \om^d \, .
$$
When we put $h=1$ we obtain 
$$
\int_{\lgl M \rgl} f*g  \cdot \om^d = \int_{\lgl M \rgl} f 
\cdot g \cdot \om^d \, .
$$
This formula means, in particular, that the functional $ \int_{\lgl M \rgl} f\cdot 
\om^d  $
is a trace functional on the deformed algebra. Such a star-products on the 
algebra of functions on a symplectic manifold were called closed star-products
in [CFS].

\vglue 1cm

\noindent {\bf Aknowledgements.} I am very much indebted to Maxim Kontsevich: 
he introduced me to the Formality conjecture in numerous informal discussions 
and he explained to me that a cyclic analog of the Formality morphism should 
exist. The work was done at the Institut des Hautes Etudes Scientifiques in 
February-March 1999, and I am grateful to the IHES for hospitality and the 
remarkable atmosphere. I am grateful 
to Mme C\'ecile Gourgues for the high quality typing 
of this manuscript.

\vglue 1cm

\noindent {\bf References}

\bigskip

\item{[CFS]} A. Connes, M. Flato, D. Sternheimer, Closed Star-Products
 and Cyclic Cohomology, Lett. in Math. Phys. 24, 1992, pp. 1-12

\smallskip 
 
\item{[K]} M. Kontsevich, Deformation quantization of Poisson manifolds, I, 
preprint q-alg/9709040.

\smallskip

\item{[L]} J.-L. Loday, Cyclic Homology, Springer-Verlag, Berlin, Die 
Grundlehren der Math. Wissenchaften, 301 (1992).

\smallskip

\item{[Sh]} B. Shoikhet, On the $A_{\ify}$-Formality conjecture, preprint
math.QA/9809117

\vglue 2cm

\noindent Boris Shoikhet

\noindent IUM

\noindent 11 Bol'shoj Vlas'evskij per.

\noindent Moscow 121002

\noindent Russia

\smallskip

\noindent e-mail address: borya@mccme.ru

\end